# Hirsch meets Fibonacci and Narayana type variants

-------------------------------------------------------------

**Leo Egghe, Hasselt University, Belgium**

leo.egghe@uhasselt.be

ORCID: 0000-0001-8419-2932

**Abstract**

For functions f of a continuous variable in $\mathbb{R}^+$ we show that the Hirsch function $h_f$ equals f iff $f(f(x)) = x\, f(x)$ on $\mathbb{R}^+$, leading for continuous f to f = **0** or the power function $f(x) = x^\alpha$, $\alpha = \frac{\sqrt{5}+1}{2}$. For functions of a discrete positive variable in $\mathbb{R}^+$, we show that $h_f = f$ implies that only the trivial function f = {(1,1)} satisfies this. We also study the problem $h_f = f \circ f$ and for $f = g \circ g$, $h_f = g$ leading to the zero function or another power law in the continuous variable case and again to f = {(1,1)} in the discrete variable case. Both problems involve the study of variants of the Fibonacci sequence for which non-trivial identities are proved and applied in the solution of the above problems.

Keywords: Fibonacci numbers; Hirsch function; power function; discrete solution; continuous solution

Mathematics Subject Classification (2020): 11B39; 11B83

## 1. Introduction

We will denote $\mathbb{R}^+ \setminus \{0\}$ as $\mathbb{R}_0^+$ and f=**0**, if the function f is equal to zero on its domain. Let now Φ be the collection of all functions of the form f:

dom(f) = $\mathbb{R}^+ \to \mathbb{R}^+$. Now $\forall \theta \in \mathbb{R}^+$, we define the Hirsch function $h_f$: $\mathbb{R}_0^+ \to \mathbb{R}^+$: $\theta \to h_f(\theta)$, where $x = h_f(\theta)$ if and only if $f(x) = \theta x$ (Egghe, 2023). Because we require $h_f$ to be a function, the equation $f(x) = \theta x$ must have a unique solution. If, however, $f(0) = 0$ then $x = 0$ is a possible second solution to this equation. This second solution is not taken into account (Egghe, 2023). The notation "h" is used because of the correspondence of the requirement $f(x) = \theta x$ with the (generalized) h-index (Egghe & Rousseau, 2019). In the case that $\theta = 1$ we obtain the (continuous version of the) famous Hirsch index (Hirsch, 2005).

This article consists of several parts. First we investigate when the continuous functions $h_f$ and f are equal. Then we consider the same problem (but with a different solution) for discrete functions. Both problems involve the use of Fibonacci numbers. Then we consider the problem when $h_f$ equals f∘f (both in the continuous and the discrete setting) Here the Narayana cows sequence is involved, featuring a characteristic determinant of this sequence. We also study when $h_{(f \circ f)}$ equals f, involving a Padovan type sequence. Relations between these sequences are proved and we show how the characteristic determinants of these sequences play a key role in the determination of variants of these sequences.

## 2. The continuous case

Considering the continuous case means that we investigate when the functions $h_f$ and f are equal on $\mathbb{R}_0^+$. We begin with the following result.

Theorem 1

If f : dom(f) = $\mathbb{R}^+ \to \mathbb{R}^+$ is continuous in 0 and in f(0) then the following expressions are equivalent:

(a) $h_f = f$ on $\mathbb{R}_0^+$ $\hspace{4cm}$ (1)

(b) $\forall x \in \mathbb{R}^+ : f(f(x)) = x f(x)$ (2)

Proof.

(a) $\Rightarrow$ (b). If $h_f = f$ on $\mathbb{R}_0^+$ then we know that for $\forall \theta \in \mathbb{R}^+$ ; f(f(θ)) = θ f(θ), which proves (b) if θ ≠ 0. By continuity we further have

$$f(f(0)) = \lim_{x \to 0^>} f(f(x)) = \lim_{x \to 0^>} x f(x) = 0 = 0 . f(0)$$

where we have moreover used that x > 0 and f(0) $\in \mathbb{R}^+$  □

(b) $\Rightarrow$ (a).  We know that $h_f(\theta)$ exists (uniquely) for every $\theta \in \mathbb{R}_0^+$. By (1) we know that $f(h_f(\theta)) = \theta\, h_f(\theta)$ (with $h_f(\theta)$ in the role of f(x)). This shows that $h_f = f$ on $\mathbb{R}_0^+$ . □

For further use, we recall the definition of Fibonacci numbers.

Definition: the Fibonacci numbers

Fibonacci numbers are whole numbers characterized by the fact that each number is the sum of the two previous ones. They are denoted as $F_k$ and exist for all k $\in \mathbb{Z}$ (the whole numbers). Choosing $F_0 = 0$ and F(1) = 1 leads to the sequence:

… -8, 5, -3, 2, -1, 1, 0, 1, 1, 2, 3, 5, 8, …..

The following results about Fibonacci numbers are well-known:

For k $\in \mathbb{N}$, $F_{-k} = (-1)^{k+1} F_k$ (3)

$\forall n \in \mathbb{Z}: F_n = \frac{\alpha^n - \beta^n}{\sqrt{5}}$ with $\alpha = \varphi = \frac{\sqrt{5}+1}{2}$ and $\beta = \frac{1-\sqrt{5}}{2}$ (4)

$\lim_{k \to +\infty} \frac{F_k}{F_{k-1}} = \varphi = \frac{\sqrt{5}+1}{2}$ (5)

$\forall n \in \mathbb{Z}: F_{n-1} F_{n+1} - F_n^2 = (-1)^n$ (6)

We further prove some lemmas and propositions in preparation for the main result of the first section.

Lemma 1. The function $h_f$ is injective on $\{\theta\,;\,h_f(\theta) \neq 0\}$

Proof. Assume that there exist $\theta_1 \neq \theta_2$ such that $h_f(\theta_1) = h_f(\theta_2) \neq 0$. As $f(h_f(\theta_1)) = \theta_1\, h_f(\theta_1)$ and $f(h_f(\theta_2)) = \theta_2\, h_f(\theta_2)$, we see that $f(h_f(\theta_1)) \neq f(h_f(\theta_2))$. This would imply that f is not a function, leading to a contradiction.

Proposition 1

If f: dom(f) = $\mathbb{R}^+ \to \mathbb{R}^+$ is continuous and satisfies (2) then either f = **0**, or f(0)=0, f(1)=1 and f is strictly increasing on $\mathbb{R}^+$.

Proof. If f satisfies (2), then, by Theorem 1, $h_f$ = f on $\mathbb{R}_0^+$ and hence by Lemma 1 f is injective on $\{\,x\,;\,f(x) \neq 0\}$. If $f(x) \neq 0$ for all x > 0, then clearly, $\{\,x\,;\,f(x) \neq 0\} = \mathbb{R}_0^+$.

By (2) we know that f(f(0)) = 0 and hence also, with f(0) in the role of x: f(f(f(0))) = f(0).f(f(0)) = 0. As f(f(0)) = 0 we see that f(0) = 0. Assume now that there exists $x_0 > 0$ such that $f(x_0) = 0$. If now $f|_{[0,x_0]} \neq \mathbf{0}$, then there exists $x_1 \in\, ]0, x_0[$ such that $f(x_1) > 0$. As f is continuous it takes, by the intermediate value theorem, all values between 0 and $f(x_1)$ on the interval [0, $x_1$] and similarly on the interval [$x_1$,$x_0$]. This shows that f is not injective on $\{\,x;\,f(x) \neq 0\}$, which is a contradiction. Hence, $f|_{[0,x_0]} = \mathbf{0}$. Now, again by the continuity of f, there exists a largest point $x_0 \in \mathbb{R}_0^+$ such that $f|_{[0,x_0]} = \mathbf{0}$ (unless, of course, f = **0**). We now consider a point $y_0 > x_0$, such that $f(y_0) > 0$. As the function f(x)/x is continuous, $f(x_0)/x_0 = 0$ and $f(y_0)/y_0 > 0$, there exists a point $y_1$ such that $0 < f(y_1)/y_1 < x_0$. Now, for $\theta = f(y_1)/y_1$ we see that $h_f(\theta) = y_1 > 0$. Indeed: $f(y_1) = \theta\, y_1$. Yet, $\theta = f(y_1)/y_1 \in [0, x_0]$. On this interval f = **0** and $h_f$ = f which contradicts $h_f(\theta) > 0$.

This shows that, unless f = **0**, f is injective on $\mathbb{R}_0^+$ and hence (Hairer & Wanner) f is strictly monotonous and thus strictly increasing. Finally, we know by (2) that f(f(1)) = f(1) and f(f(f(1))) = f(1).f(f(1)). Applying the first equality in the second one leads to f(f(1)) = f(1).f(f(1)). As f(f(1)) > 0 (f is strictly increasing), it follows that f(1) = 1. □

The next proposition is the crux of the proof of our most important result of this section (Theorem 2) and makes use of Fibonacci's numbers. We already note that intuitively one might expect a relationship with the Fibonacci numbers as in the equality f(f(x))= f$^{(2)}$(x),  f$^{(2)}$(x) is the product of its two preceding ones, namely f$^{(0)}$(x) = x and f$^{(1)}$(x)=f(x).

Proposition 2

If f : dom(f) = $\mathbb{R}^+ \to \mathbb{R}^+$ is continuous then (2) implies the following statement (7)

$$\forall x \in \mathbb{R}^+, \forall\, n \in \mathbb{Z}, \text{and with } y = f(x):$$

$$f(x^{F_n} y^{F_{n+1}}) = x^{F_{n+1}} y^{F_{n+2}} \tag{7}$$

Proof (based on Edgar's answer in https://math.stackexchange.com/questions/1412525/find-all-functions-f-such-that-ffx-fxx/1412759#1412759 ).  We note that by Proposition1, y= f(x) > 0. For fixed n we refer to $f(x^{F_n} y^{F_{n+1}}) = x^{F_{n+1}} y^{F_{n+2}}$ as (ρ$_n$). Then f(x) = y is nothing but (ρ$_{-1}$). By (2) we know that f(f(x)) = x f(x), hence f(y) = xy, which is nothing but (ρ$_0$). Next, we prove (ρ$_1$). We know, by (2) that f(f(f(x)))= f(x)f(f(x)). Now the left-hand side is equal to f(xf(x)) = f(xy), while the right-hand side is yf(y), which is xy$^2$. This shows (ρ$_1$), namely f(xy) = xy$^2$. We will not write a formal induction step (for natural numbers) but just go one step further to illustrate the method.

We know that f(f(f(f(x)))) = f(f(x)).f(f(f(x))) which can be rewritten as $f(xy^2)$ = (xy)(xy²), proving (ρ₂), namely: $f(xy^2) = x^2y^3$.

So, we have shown (7) for all natural numbers and for n = -1. Now we consider the negative whole numbers n = -k, k ∈ ℕ.

Using (3), equation (7) leads to the following equalities that we have to prove for all x ∈ $\mathbb{R}_0^+$ and y= f(x):

$$f(x^{F_{-k}} y^{F_{-k+1}}) = x^{F_{-k+1}} y^{F_{-k+2}}$$

$$\Leftrightarrow f\left(x^{(-1)^{k+1}F_k} y^{(-1)^k F_{k-1}}\right) = x^{(-1)^k F_{k-1}} y^{(-1)^{k-1} F_{k-2}}$$

$$\Leftrightarrow f\left(\frac{y^{(-1)^k F_{k-1}}}{x^{(-1)^k F_k}}\right) = \frac{x^{(-1)^k F_{k-1}}}{y^{(-1)^k F_{k-2}}} \tag{8}$$

We know from the definition of the generalized h-index that $h_f(\theta) = x$ iff $\theta = \frac{f(x)}{x}$, assuming that the generalized h-index exists (uniquely) for this value of θ. This implies f(y/x)=x. Another argument to show the same equality runs as follows: by (2) and Proposition 1, we know that $h_f = f$, f is strictly increasing and hence $f^{-1}$ is well-defined. Applying (2) to $f^{-1}(x)$ in the role of x leads to y = f(x) = x. $f^{-1}(x)$ or $f^{-1}(x)$ = y/x or x = f(y/x), which is (ρ₋₂).

Now we apply (2) to $f^{-1}(f^{-1}(x)) \in \mathbb{R}^+$, leading to

x = $f^{-1}(x).f^{-1}(f^{-1}(x))$ which, using (ρ₋₂), gives: x = (y/x). $f^{-1}(f^{-1}(x))$, or

$x^2/y = f^{-1}(f^{-1}(x))$, hence $f(x^2/y) = f^{-1}(x)$ = y/x and, finally

$$f\left(\frac{x^2}{y}\right) = \frac{y}{x} \tag{ρ₋₃}$$

As, before, we do not write down the formal induction step, being convinced that this is clear to the reader. This proves (8), hence (7) for all integers (ℤ). □

Reflecting on (8) we see that it shows that f transfers y in the numerator and x in the denominator, to x in the numerator and y in the denominator, and vice versa. For instance, take n = -k ∈ $\mathbb{Z}^-$, with k even leading to

$$f\left(\frac{y^{F_{k-1}}}{x^{F_k}}\right) = \frac{x^{F_{k-1}}}{y^{F_{k-2}}} \tag{9}$$

with a similar expression in the case of k odd.

Finally, we formulate and prove the main theorem (Theorem 2) of this first section.

Theorem 2

If f: dom(f) = $\mathbb{R}^+$ → $\mathbb{R}^+$ is continuous then the following two expressions are equivalent:

(a) $\forall\, x \in \mathbb{R}^+$: f(f(x)) = x f(x)  (10)

(b) f = **0** (the null function) or $f(x) = x^\varphi$ with $\varphi = \frac{\sqrt{5}+1}{2}$  (11)

Here we recognize φ as the golden ratio.

Proof. We first prove that (b) implies (a).

If f=**0** then (a) is trivially satisfied. We next consider f(x) = $x^p$ ($x \in \mathbb{R}^+$ and $p \in \mathbb{R}$ ). If (2) holds for f(x) = $x^p$ then $x^{(p^2)} = x\, x^p$. As this equality holds for all x ∈ $\mathbb{R}^+$, it follows that $p^2-p-1=0$ and thus, $p = \frac{1\pm\sqrt{5}}{2}$. The case $p = \frac{1-\sqrt{5}}{2}$ < 0 cannot occur as then f(0) ≠ 0, which is required by Proposition 1. Hence, a function $x^p$ with f(f(x)) = xf(x) exists only if p = $\varphi = \frac{\sqrt{5}+1}{2}$.

Proof of the implication (a) ⟹(b).

We first rewrite (5) as $\lim_{k \to \infty}(\varphi F_{k-1} - F_k) = 0$. If now $y = f(x) > x^\varphi$ for some x then $\frac{y^{F_{k-1}}}{x^{F_k}} > x^{(\varphi F_{k-1} - F_k)}$, which tends to $x^0 = 1$ by (5).

Further, $\frac{x^{F_{k-1}}}{y^{F_{k-2}}} < x^{(F_{k-1} - \varphi F_{k-2})}$, which also tends to $x^0 = 1$ by (5). This shows that, by (9) (or (7)), f transforms values strictly larger than 1 to values strictly smaller than 1 (consider (8)), unless

$$\lim_{\substack{k \to \infty \\ k \text{ even}}} \left(\frac{y^{F_{k-1}}}{x^{F_k}}\right) = 1 = \lim_{\substack{k \to \infty \\ k \text{ even}}} \left(\frac{x^{F_{k-1}}}{y^{F_{k-2}}}\right)$$

Let now $a > 1$ such that $y = ax^\varphi > x^\varphi$. Then

$$\frac{y^{F_{k-1}}}{x^{F_k}} = a^{F_{k-1}} x^{(\varphi F_{k-1} - F_k)}$$

The first factor in this expression tends to infinity, for k tending to infinity, while the second one tends to 1 for k (even) tending to infinity. Hence $\lim_{\substack{k \to \infty \\ k \text{ even}}} \left(\frac{y^{F_{k-1}}}{x^{F_k}}\right) = +\infty$. Similarly in $\left(\frac{x^{F_{k-1}}}{y^{F_{k-2}}}\right) = \frac{1}{a^{F_{k-2}}} x^{(F_{k-1} - \varphi F_{k-2})}$ the first factor tends to zero when k tends to infinity and the second one tends to 1 (with k even). Hence: $\lim_{\substack{k \to \infty \\ k \text{ even}}} \left(\frac{x^{F_{k-1}}}{y^{F_{k-2}}}\right) = 0$.

By choosing an odd value for k we can similarly prove that $y = f(x) < x^\varphi$ is not possible. This shows that $\forall x \in \mathbf{R}^+: f(x) = x^\varphi$ with $\varphi = \frac{\sqrt{5}+1}{2}$, the golden ratio. □

## 3 The discrete case

A general discrete function can be denoted as f: $\{a_1, \ldots, a_n\} \to \{b_1, \ldots b_m\}$, where we assume that all $a_j > 0$. We have that m ≤ n otherwise f is not a function. If m < n then f is not injective and hence $h_f$ (= f) is not injective, which is not possible by Lemma 1 (also valid here), and since all $a_j \neq 0$, being the images of $h_f$ (see also below). Hence we have f: $\{a_1, ...,a_n\} \rightarrow \{b_1,...b_n\}$ with $f(a_j) = b_j$, for $j=1,..., n$. Such a function can be represented as a scatterplot, see Fig.1.

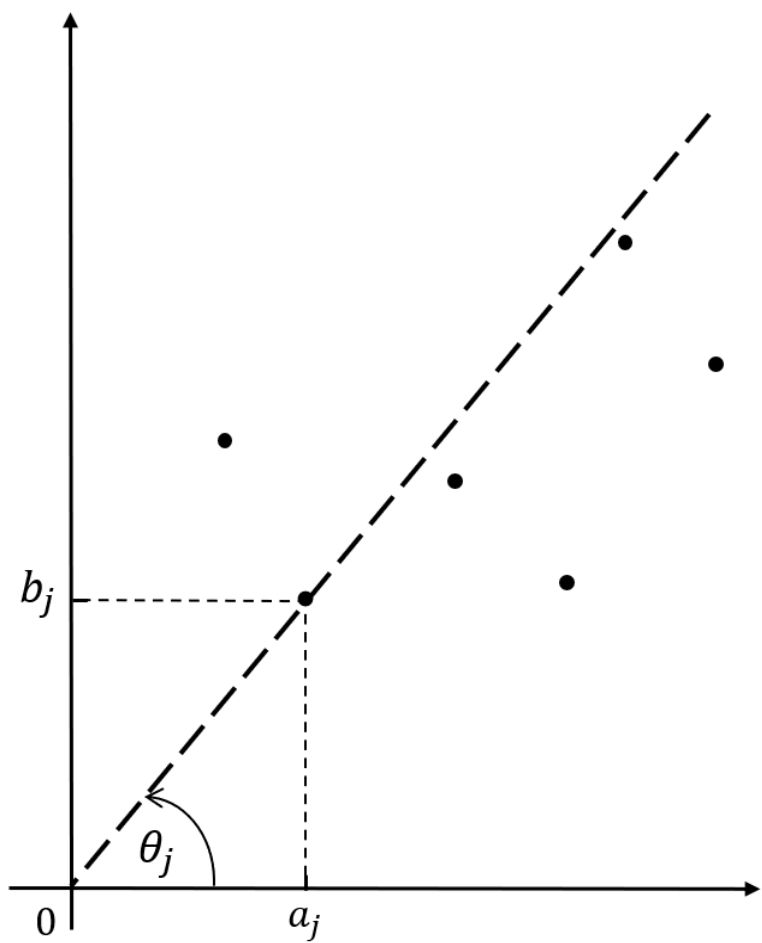

Fig. 1. Discrete function f

Now the function $h_f(\theta)$ is defined for values $\theta_j = f(a_j)/a_j = b_j/a_j$ $j = 1,...,n$, (and only these) and maps $\theta_j$ to $a_j$. Now the requirement $h_f = f$ can only be valid if both functions have the same domain. As $h_f$ is defined on values of the form $f(a_j)/a_j$ and the domain of f consists of all a-values, the equality between $h_f$ and f can only hold if, for all $a_i \in \{a_1, ...,a_n\}$, there exists j such that $a_i = f(a_j)/a_j$ and then $h_f(f(a_j)/a_j) = a_j = f(a_i) = b_i$. This shows that each b-value is equal to some a-value, hence $\{b_1,…b_n\} \subset \{a_1, …,a_n\}$. As $h_f$ is defined in points of the form $f(a_j)/a_j$ and takes values $a_j$ it follows that $h_f$ is surjective (all a-values are different), leading to $\{b_1,…b_n\} = \{a_1, …,a_n\}$, from which it follows that f is a permutation of the set $\{a_1, …,a_n\}$.

In the proof of the next theorem, we will need the following lemma about Fibonacci numbers.

Lemma 2

$$\forall n > 2: (F_{n-1} - 1)(F_{n+1} - 1) < F_n^2 \qquad (12)$$

Proof. We first note that for n > 2, $F_{n-1} + F_{n+1} > (-1)^n + 1$. Now

$(F_{n-1} -1)(F_{n+1} -1) = F_{n-1} F_{n+1} - F_{n-1} - F_{n+1} + 1 =$ (by (6)) $(-1)^n + F_n^2 - F_{n-1} - F_{n+1} + 1 < (-1)^n + F_n^2 + (-1)^{n+1} - 1 + 1 = F_n^2$ □

We next prove the following theorem.

Theorem 3

$h_f = f$ if and only if n =1 and f maps 1 to 1.

Proof. If f has fixed points, then these are removed from the set $\{a_1, …,a_n\}$, leading to a function f* for which $h_{f^*} = f^*$ holds. We now denote f* again as f and prove now that f is the empty function (ø). Then the original function f has only fixed points. But from $a_i = f(a_j)/a_j$, it follows that $a_i = 1$. Hence n =1 and f maps 1 to 1.

First, we assume that $f(a_1) = a_2$; $f(a_2) = a_3$; …, $f(a_{n-1}) = f(a_n)$ and $f(a_n) = a_1$.

The relation $h_f = f$ implies $f(f(a_j)) = a_j.f(a_j)$ for every j = 1, …, n. We first show, as an illustration, that f = ø for the simple cases n=1 and n=2. For n=1 we see that f = ø because f has no fixed point. For n=2 and assuming that f ≠ ø, we have $f(a_1) = a_2$ and $f(a_2) = a_1$. But, we also have from $f(f(a_1)) = a_1.f(a_1)$ that $f(a_2) = a_1 = a_1.a_2$. Hence $a_2 = 1$, and similarly $a_2 = a_2.a_1$ and hence $a_1 = 1$. This, however, is not possible as f has no fixed points. We conclude that f = ø.

Now we consider the general case n > 2. For the permutation f and using $f(f(a_j)) = a_j.f(a_j)$ for every j = 1, …, n, we have:

$$a_3 = a_1a_2,\ a_4=a_2a_3,\ …,\ a_n = a_{n-1}a_{n-2},\ a_1=a_{n-1}a_n,\ a_2 = a_na_1. \tag{13}$$

We will show that from these equalities it follows that all $a_j$, j=1, …n must be equal to 1 (see further), which is excluded by the requirement that f must be a function without fixed points. Hence f = ø.

It follows from (13), read from left to right, that

$a_3 = a_1a_2$ or $a_3 = a_1^{F_1}a_2^{F_2}$; $a_4=a_1a_2^2$ or $a_4 = a_1^{F_2}a_2^{F_3}$; $a_5 = a_1^2a_2^3$ or $a_5 = a_1^{F_3}a_2^{F_4}$; $a_6 = a_1^3a_2^5$ or $a_6 = a_1^{F_4}a_2^{F_5}$; …; $a_k = a_1^{F_{k-2}}a_2^{F_{k-1}}$ (k=3,…, n) and finally $a_1 = a_1^{F_{n-1}}a_2^{F_n}$ and $a_2 = a_1^{F_n}a_2^{F_{n+1}}$.

This leads to: $a_1^{1-F_{n-1}} = a_2^{F_n}$ and $a_2^{1-F_{n+1}} = a_1^{F_n}$, hence $a_2 = a_1^{\left(\frac{F_n}{1-F_{n+1}}\right)}$, and thus $a_1^{1-F_{n-1}} = a_1^{\left(\frac{F_n^2}{1-F_{n+1}}\right)}$. As by (12) $(F_{n-1} - 1)(F_{n+1} - 1) < F_n^2$, it follows that $a_1 = 1$. It then follows from $a_2 = a_1^{F_n}a_2^{F_{n+1}}$ that also $a_2 = 1$ and then from (13) that all $a_j$, j=1, …n, are equal to 1. As f must be a function without fixed points, it follows that f = ø.

Finally, the above proof of the result for the case that $f(a_1) = a_2$; $f(a_2) = a_3$; …, $f(a_{n-1}) = f(a_n)$ and $f(a_n) = a_1$ suffices. Indeed, if f is any permutation without fixed points, and by changing the notation, we may assume that $f(a_1) = a_2 \neq a_1$. We know that $f(a_2) \neq a_2$ (no fixed points). Assume now that $f(a_2) = a_1$ then we can apply the previous result for n=2 and conclude that $f|_{\{a_1,a_2\}}$= ø and hence f itself is equal to ø. This application is possible since, for any loop A of f (as {$a_1, a_2$} here) we have $h_{f|_A} = (h_f)|_A$ which follows from the definition of $h_f$ and since A is a loop: $x = (h_{f|_A})(\theta) \Leftrightarrow (f|_A)(x) = \theta x \Leftrightarrow f(x) = \theta x$ and $x \in A$. Now $\theta \in A$ since $h_f(\theta) = f(\theta) = x \in A$, and A is a loop. Otherwise, assume that $f(a_2)$ is equal to another $a_j$, which we rename $a_3$. If $f(a_3) = a_1$ we have a cycle of order three and applying the previous result for n=3, again leads to f = ø. This continues until we either find that f = ø or f is exactly the permutation of the previous result and hence we always have f = ø.

**4. A variant of the Fibonacci numbers: the Narayana cows sequence $(E_n)_{n \in \mathbb{Z}}$**

Instead of the requirement $h_f = f$, we now study the requirement $h_f = f \circ f$. Similar to Theorem 1, we prove Theorem 4.

Theorem 4

If f : dom(f) = $\mathbb{R}^+ \to \mathbb{R}^+$ is continuous then the following expressions are equivalent:

(a) $h_f = f \circ f$ on $\mathbb{R}_0^+$ $\hspace{4cm}$ (14)

(b) $\forall\, x \in \mathbb{R}^+ : f\big(f(f(x))\big) = x\, f(f(x))$ $\hspace{2cm}$ (15)

Proof.

(a) $\Rightarrow$ (b). If $h_f = f \circ f$ on $\mathbb{R}_0^+$ then we know that for $\forall\, \theta \in \mathbb{R}^+$, $f(f(f(\theta))) = \theta\, f(f(\theta))$, which proves (b) if $\theta \neq 0$. By continuity we further have

$$f\big(f(f(0))\big) = \lim_{\substack{x \to 0 \\ >}} f\big(f(f(x))\big) = \lim_{\substack{x \to 0 \\ >}} x\, f(f(x)) = 0 = 0 \cdot f(f(0))$$

(b) ⇒(a). We know that $h_f(\theta)$ exists (uniquely) for every $\theta \in \mathbb{R}_0^+$ and that it is characterized by $f(h_f((\theta)) = \theta\, h_f(\theta)$. Hence, it follows from (15) that $h_f$ = f∘f on $\mathbb{R}_0^+$. □

Theorem 5. Expression (15) implies that f is either strictly increasing with f(0) = 0, f(1) = 1, or f = **0**.

Proof. From (15) we see that f(f(f(0))) = 0 (take x =0) and f(f(f(f(0)))) = f(0).f(f(f(0))) (take x = f(0)). Then, the left-hand side is equal to f(0), while the right-hand side is equal to 0. Hence f(0) = 0. Assume now that there exists x > 0 such that f(x) = 0. Then f(f(0)) = 0 and f(f(x)) = 0. Yet, as $h_f$ = f∘f we know that the function f∘f is injective on the set { f∘f ≠ 0}. Consequently, (f∘f)|$_{[0,x]}$ = **0** (on this interval). If f ≠ **0**, then there exists a largest strictly positive point x such that (f∘f)|$_{[0,x]}$ = **0**, and thus we have for all y >x: (f∘f) > 0. Then, as f(0) = 0, also (f∘f)(y) > 0. For y > x such that $0 < \frac{f(y)}{y} = \theta \leq x$ and hence $h_f(\theta)$ = y > 0, which is impossible as θ ∈ [0,x], $h_f$ = f∘f and (f∘f)|$_{[0,x]}$ = **0**. We hence see that, unless f = **0**, f∘f is injective on $\mathbb{R}_0^+$; then also f in injective on $\mathbb{R}_0^+$ and thus strictly monotone. Then it follows from f(0) = 0 that f is strictly increasing. Finally, (15) implies that f(f(f(1))) = f(f(1)) (by taking x = 1) and f(f(f(1))) = f(1) f(f(f(1))) (by taking x = f(1)). The left-hand side is equal to f(f(f(1))), which implies that f(1) =1 unless f(f(1)) = 0, but this is impossible as f(0) = 0 and the function f is strictly increasing. □

Proposition 3. $\forall\, x \in \mathbb{R}^+ : f\left(f(f(x))\right) = x\, f(f))$, i.e., expression (15), is satisfied for f = **0** and for $f(x) = x^{a_0}$, with the number $a_0$ being the real-valued solution of $s^3 - s^2 - 1 = 0$.

Proof. It is trivial to see that if a is a solution of $s^3 - s^2 - 1 = 0$, then $x^{a_0}$ satisfies equation (15).

We note that $s^3 - s^2 - 1 = 0$ has the following (rounded) solutions:

$a_0 \approx 1.4648493$ (the real-valued solution) and

$a_{2,3} \approx -0.2324247 \pm i(0.7931913) = \rho(\cos(\zeta) \pm i \sin(\zeta))$ with $\rho \approx 0.8265432$ and $\zeta \approx 1.8558416$ (radians).

Theorem 6

$\forall\, x \in \mathbb{R}^+ : f\big(f(f(x))\big) = x\, f(f(x))$, i.e., expression (15) implies that $\forall\, x \in \mathbb{R}^+, y = f(x), \forall\, n \in \mathbb{Z}$:

$$f(x^{E_n} y^{E_{n-1}} f(y)^{E_{n+1}}) = x^{E_{n+1}} y^{E_n} f(y)^{E_{n+2}} \qquad (\varepsilon_n)$$

where $(E_n)_n$ is the unique sequence satisfying the relation

$$\forall\, n \in \mathbb{Z}: E_n + E_{n+2} = E_{n+3} \qquad (16)$$

with $E_0 = 0$, $E_1 = E_2 = 1$.

Before proving this theorem, we show Table 1 which contains the numbers $E_n$ for n = -16 to 13

Table 1. E-numbers

| n | -1 | 0 | 1 | 2 | 3 | 4 | 5 | 6 | 7 | 8 | 9 | 10 | 11 | 12 | 13 |
|---|---|---|---|---|---|---|---|---|---|---|---|---|---|---|---|
| $E_n$ | 0 | 0 | 1 | 1 | 1 | 2 | 3 | 4 | 6 | 9 | 13 | 19 | 28 | 41 | 60 |
| n | -16 | -15 | -14 | -13 | -12 | -11 | -10 | -9 | -8 | -7 | -6 | -5 | -4 | -3 | -2 |
| $E_n$ | 9 | 4 | -8 | 1 | 5 | -3 | -2 | 3 | 0 | -2 | 1 | 1 | -1 | 0 | 1 |

Proof (of Theorem 6)

We note that f(y) = f(y) is ($\varepsilon_{-1}$) and f(x) = y is ($\varepsilon_{-2}$).

We first prove that (15) implies ($\varepsilon_0$). Indeed, it follows immediately that f(f(y)) = x.f(y), which is ($\varepsilon_0$).

Next, we replace in expression (15) x by f(x) leading to: f(f(f(f(x)))) = f(x).f(f(f(x))) or f(xf(y))=y.xf(f(x)) = yx f(y), which exactly is ($\varepsilon_1$).

Now, replacing x in (15) by f(f(x)) leads to:

f(f(f(f(f(x))))) = f(f(x)).f(f(f(f(x)))), using ($\varepsilon_1$) this becomes f(xyf(y)) = f(y). xyf(y) and hence f(xyf(y)) = xy (f(y))$^2$ i.e., equality ($\varepsilon_2$).

Now we replace x in (15) by f$^{-1}$(x), leading to:

f(f(f(f$^{-1}$(x))))= f$^{-1}$(x). f(f(f$^{-1}$(x))) or f(y) = f$^{-1}$(x). y. Hence f$^{-1}$(x) = f(y)/y. Applying f to this expression leads to x = f$^{-1}$(y$^{-1}$.f(y)) which is ($\varepsilon_{-3}$).

It is now clear how to continue using (16). □

We next mention the analog of Binet's formula (4) for the sequence $(E_n)_n$ It suffices to solve (16) which is a homogeneous linear difference equation with constant coefficients and initial conditions $E_0 = 0$, $E_1 = E_2 = 1$.

For $n \in \mathbb{Z}$ we find:

$$E_n = C_1 a_0^n + C_2 \rho^n \cos(n\zeta) + C_3 \rho^n \sin(n\zeta) \tag{17}$$

with $a_0$, $\rho$, and $\zeta$ as in proposition 3. Inserting the initial conditions leads to $C_1$ = 0.4173489 = -$C_2$ and $C_3$ = 0.367685.

The next theorem yields limit properties, which will be used further on.

Theorem 7

(i) $\lim_{n\to+\infty}(a_0 E_n - E_{n+1}) = 0$ and hence (i') $\lim_{n\to+\infty}\frac{E_{n+1}}{E_n} = a_0$ (18)

(ii) $\lim_{n\to-\infty}(E_n + a_0 E_{n-1} + a_0^2 E_{n+1}) = 0$ (19)

Proof

For all natural numbers n we have:

$$a_0 E_n - E_{n+1} = C_1 a_0^n - C_1 a_0 \rho^n \cos(n\zeta) + C_3 a_0 \rho^n \sin(n\zeta)$$
$$- C_1 a_0^n + C_1 \rho^{n+1} \cos((n+1)\zeta) - C_3 \rho^{n+1} \sin((n+1)\zeta)$$

As 0 < ρ < 1 and |sin(x)| ≤ 1 as well as | cos(x)| ≤ 1, the limits in (18) follow easily.

Now, from Proposition 3 and Theorem 6, namely that $f(x) = x^{a_0}$ satisfies ($\varepsilon_k$) for all k, we write, using ($\varepsilon_{k-1}$):

$$x^{a_0(E_{k-1} + a_0 E_{k-2} + a_0^2 E_k)} = x^{(E_k + a_0 E_{k-1} + a_0^2 E_{k+1})}$$

Now we keep k fixed and write this equality for the cases ($\varepsilon_{k-2}$), …, ($\varepsilon_{k-m}$), with m taking increasing positive values. This leads to:

$$a_0^m (E_{k-m} + a_0 E_{k-m-1} + a_0^2 E_{k-m+1}) = E_k + a_0 E_{k-1} + a_0^2 E_{k+1}$$

and hence:

$$E_{k-m} + a_0 E_{k-m-1} + a_0^2 E_{k-m+1} = \frac{1}{a_0^m}(E_k + a_0 E_{k-1} + a_0^2 E_{k+1})$$

We note that the last factor is fixed, while $\frac{1}{a_0^m}$ tends to zero when m tends to +∞. When m tends to +∞, then k-m tends to -∞ (k is fixed). Finally, writing k-m as n yields (19). □

Now we come to the crux of our theory.

Theorem 8

The expressions $f(x^{E_n} y^{E_{n-1}} f(y)^{E_{n+1}}) = x^{E_{n+1}} y^{E_n} f(y)^{E_{n+2}}$ ($\varepsilon_n$)

from Theorem 6 can only occur for f = **0** and for f(x) = $x^{a_1}$ with $a_0$ = 1.4648493…, i.e., the only real solution of $s^3 - s^2 - 1 = 0$. We will show that for f(x) ≠ $x^{a_0}$ and not f = **0**, for |n| high, n ∈ $\mathbb{Z}^-$, ($\varepsilon_n$) implies that f transforms numbers strictly smaller than one into numbers strictly larger than one and vice versa. This, however, is an impossibility as f is strictly increasing with f(1) = 1 (by Theorem 5).

We note that for |n| high $a_1^n \approx 0$ hence sgn($E_n$) is given by:

$$sgn(-C_1\cos(n\zeta) + C_3\sin(n\zeta)) = sgn(-C_1\cos(m\zeta) - C_3\sin(m\zeta))$$

with m = |n| = -n. We next investigate the changes in the sign of $c_m$, with

$$c_m = \frac{C_1}{C_3}\cos(m\zeta) + \sin(m\zeta) \qquad (20)$$

Table 2 shows the values of $c_m$ for m =1, …, 17 (we do not need the values for higher m). Recall that the values for $C_1$, $C_3$, and $\zeta$ are known.

Table 2. Approximate values and signs of $c_m$

| m | $c_m$ | m | $c_m$ |
|---|---|---|---|
| 1 | 0.6404657 | 9 | -1.4569673 |
| 2 | -1.495271 | 10 | 0.8002606 |
| 3 | 0.2004774 | 11 | 1.0068992 |
| 4 | 1.3825221 | 12 | -1.3665426 |
| 5 | -0.9780103 | 13 | -0.2383533 |
| 6 | -0.8324874 | 14 | 1.5005929 |
| 7 | 1.4462027 | 15 | -0.6055829 |
| 8 | 0.0191404 | 16 | -1.160012 |
|  |  | 17 | 1.1259757 |

The change in signs illustrated in Table 2 is repeated indefinitely. Indeed, consider one value of m $\zeta$ (m a strictly positive natural number). We search a strictly positive natural number m' such that (m' $\zeta$) is arbitrarily close to (m ζ) + 2kπ (k ∈ ℕ), i.e., (m' ζ) ≈ (m ζ) + 2kπ or m' ≈ m + 2kπ/ ζ. We know that there exist rational numbers of the form c/d (c,d ∈ ℕ) that approximate π/ ζ arbitrarily close. Hence, we can take m'= m + 2kc/d for every k which is a multiple of d. This leads to infinitely many m' for which the change in signs as shown in Table 1, is repeated. From Table 2 we see (for m = 9 = -n, and hence for infinitely many n) the following changes in sign

| $E_{n-1}$ | $E_n$ | $E_{n+1}$ |
|---|---|---|
| + | - | + |
| $E_{n+4}$ | $E_{n+5}$ | $E_{n+6}$ |
| - | + | - |

Now $(\varepsilon_n)$ implies that

$$\underbrace{(f \circ \ldots \circ f)}_{5\ times}(x^{E_n} y^{E_{n-1}} f(y)^{E_{n+1}}) = x^{E_{n+5}} y^{E_{n+4}} f(y)^{E_{n+6}} \qquad (21)$$

Hence, if y = f(x) > $x^{a_0}$ and f(y) > $x^{(a_0^2)}$ (with x fixed) we have, that

$x^{E_n} y^{E_{n-1}} f(y)^{E_{n+1}} > x^{(E_n + a_0 E_{n-1} + a_0^2 E_{n+1})}$ which tends to $x^0$ = 1 by (19). Moreover, we have, again by (19), that

$x^{E_{n+5}} y^{E_{n+4}} f(y)^{E_{n+6}} < x^{(E_{n+5} + a_0 E_{n+4} + a_0^2 E_{n+6})}$ which also tends to $x^0$ = 1. This shows that the function $\underbrace{(f \circ \ldots \circ f)}_{5\ times}$ maps numbers larger than one to numbers smaller than one, which is not possible by Theorem 5, unless

$\lim_{n \to -\infty}(x^{E_n} y^{E_{n-1}} f(y)^{E_{n+1}}) = \lim_{n \to -\infty}(x^{E_{n+5}} y^{E_{n+4}} f(y)^{E_{n+6}})$ = 1. Yet, this equality is impossible. Indeed, from y = f(x) > $x^{a_0}$ and f(y) > $x^{(a_0^2)}$ we can write: f(x) = r $x^{a_0}$ and f(y) = s $x^{(a_0^2)}$, with r,s > 1 fixed (as x is fixed). Then:

$$\lim_{n \to -\infty}(x^{E_n} y^{E_{n-1}} f(y)^{E_{n+1}}) = \lim_{n \to -\infty}\left(r^{E_{n-1}} s^{E_{n+1}} x^{(E_n + a_0 E_{n-1} + a_0^2 E_{n+1})}\right)$$

The product of the first two factors is strictly larger than one, while the limit for n tending to -∞ of the last factor is equal to one, (by (19)). Hence this limit is strictly larger than one. Moreover, we have:

$$\lim_{n \to -\infty} (x^{E_{n+5}} y^{E_{n+4}} f(y)^{E_{n+6}}) = \lim_{n \to -\infty} \left( r^{E_{n+4}} s^{E_{n+6}} x^{(E_{n+5}+a_0 E_{n+4}+a_0^2 E_{n+6})} \right)$$

As $E_{n+4}$ and $E_{n+6}$ are strictly smaller than zero and are integers, it follows that they are smaller than or equal to -1. Then the product of the first two factors is smaller than or equal to 1/(rs) and thus strictly smaller than one. The limit for n tending to -∞ of the last factor is equal to one, (by (19)). Hence this limit is strictly smaller than one. This excludes the case y = f(x) > $x^{a_0}$ and f(y) > $x^{(a_0^2)}$. Similarly we can exclude the case y = f(x) < $x^{a_0}$ and f(y) < $x^{(a_0^2)}$. Assume now that y = f(x) > $x^{a_0}$ and f(y) ≤ $x^{(a_0^2)}$. We now need Table 2 for m = 5 = -n (and hence for infinitely many n), leading to:

| $E_{n-1}$ | $E_n$ | $E_{n+1}$ |
|---|---|---|
| + | - | - |
| $E_{n+1}$ | $E_{n+2}$ | $E_{n+3}$ |
| - | + | + |

Applying the same argument as above and using (22) instead of (21)

$$(f \circ f)( (x^{E_n} y^{E_{n-1}} f(y)^{E_{n+1}}) = x^{E_{n+2}} y^{E_{n+1}} f(y)^{E_{n+3}} \qquad (22)$$

where (22) results from applying ($\varepsilon_n$) twice, leads to the conclusion that y = f(x) > $x^{a_0}$ and f(y) ≤ $x^{(a_0^2)}$ is impossible. Finally, we also have (same argument) that y = f(x) < $x^{a_0}$ and f(y) ≥ $x^{(a_0^2)}$ is impossible. Hence f(x) < $x^{a_0}$ and f(x) > $x^{a_0}$ are both impossible. Hence, we conclude that for x ∈ $\mathbb{R}^+$, f(x) = $x^{a_0}$. □

## 5. Generalizations of the above theory

5.1 Using more than two compositions

A simple idea to generalize the above theory is to use more than two compositions leading to:

$$h_f = \underbrace{f \circ f \circ \ldots \circ f}_{n \text{ times}} \text{ on } \mathbb{R}_0^+ \tag{23}$$

with $n \in \mathbb{N}_0$. The cases n=1 and n=2 are already studied. As in Theorem 4 we have:

$$(\underbrace{f \circ f \circ \ldots \circ f}_{(n+1) \text{ times}})(x) = x(\underbrace{f \circ f \circ \ldots \circ f}_{n \text{ times}})(x) \tag{24}$$

The functions f = **0** and f=$x^c$ are solutions with

$$c^{n+1} - c^n - 1 = 0 \tag{25}$$

5.2 Another interesting, related problem is the following:

Given a continuous function g: $\mathbb{R}^+ \to \mathbb{R}^+$ find f such that $f = \underbrace{g \circ g \circ \ldots \circ g}_{m \text{ times}}$, m $\in \mathbb{N}$, and

$$h_f = h_{(\underbrace{g \circ g \circ \ldots \circ g}_{m \text{ times}})} = g \tag{26}$$

Definition

If $f = \underbrace{g \circ g \circ \ldots \circ g}_{m \text{ times}}$ then we say that for n = 1/m, g is the n-fold composition of f with itself, denoted as $g =: f^{(1/n)} = \underbrace{f \circ \ldots \circ f}_{n \text{ times}}$. In this way, we have defined $g = f^{(1/n)} = \underbrace{f \circ \ldots \circ f}_{n \text{ times}}$, for all $n \in \mathbb{Z} \setminus \{0\}$.

Using this definition, the problem in (26) is: find f such that

$$h_f = \underbrace{f \circ \ldots \circ f}_{n \text{ times}} \tag{27}$$

with n = 1/m, m $\in \mathbb{N}$. This problem generalizes (23) and leads to (cf. (26))

$$(\underbrace{g \circ g \circ \ldots \circ g}_{m+1 \ times})(x) = x \, g(x) \qquad (28)$$

Again $g = 0$ and $g = x^c$ satisfy (28) for c the real solution of $c^{m+1} - c - 1 = 0$. This yields the solutions $f = \mathbf{0}$ and $f(x) = (\underbrace{g \circ g \circ \ldots \circ g}_{m \ times})(x) = x^{(c^m)}$.

5.3 The special case n = -m = -2.

In this special case, we have to find $f = g \circ g$ such that on $\mathbb{R}_0^+$,

$$h_{(g \circ g)} = g \qquad (29)$$

or, for all x in $\mathbb{R}^+$:

$$g(g(g(x))) = x \, g(x) \quad \text{and} \quad f = g \circ g \qquad (30)$$

The equations (30) are satisfied for $g = \mathbf{0}$ and $g(x) = x^c$, where c is the real solution of

$$c^3 - c - 1 = 0 \qquad (31)$$

and $f(x) = x^{(c^2)}$. We find one real solution for (30) and two complex conjugate ones:

$$c = c_1 \approx 1.3247178 \qquad (32)$$

$c_{2,3} \approx -0.6623589 \pm i \, (0.5622796)$

This ends the analog of Proposition 3. Next, we consider the analog of Theorem 6.

Theorem 9

The equality $g(g(g(x))) = x \, g(x)$ implies that $\forall \, x \in \mathbb{R}^+, y = g(x), \forall \, n \in \mathbb{Z}$:

$$g(x^{D_n} y^{D_{n+2}} g(y)^{D_{n+1}}) = x^{D_{n+1}} y^{D_{n+3}} g(y)^{D_{n+2}} \qquad (\delta_{n-2})$$

with $(D_n)_n \in \mathbb{Z}$ the unique sequence determined by

$$D_n + D_{n+1} = D_{n+3} \qquad (33)$$

with $D_{-1}=-1$, $D_0=1$, $D_1=D_2=0$.

We have then further that

$$D_n = C_1' c_1^n + C_2'(\rho')^n \cos(n\zeta') + C_3'(\rho')^n \sin(n\zeta') \qquad (34)$$

with $c_1 \approx 1.3247178 > 1$. This leads to the analog of Theorem 7 (without proof). From these results, it follows that the equalities $(\delta_n)$, $n \in \mathbb{Z}$, can only occur for g = **0** or $g(x) = x^{c_1}$ and thus (28) only occurs for f = g∘g=**0** or $f(x) = (g \circ g)(x) = x^{(c_1^2)} = x^{1.7548772}$ (approximately).

We end this section with the following conjecture.

We conjecture that on $\mathbb{R}^+$ (27) only has the solutions f=0 and $f(x) = x^{d_n}$, with $d_n > 1$ decreasing in n. Recall that we already know that $d_1 = \frac{\sqrt{5}+1}{2} \approx$ 1.618034 ; $d_2 = a_0 \approx 1.4648496$ and $d_{(1/2)} = c^2 \approx 1.7548772$.

## 6. The problem $h_f = f \circ f$ in the discrete case

Similar to the case $h_f = f$ we may assume that f maps $\{a_1, \ldots, a_n\}$ to itself with

$$f(a_1) = a_2, f(a_2) = a_3, \ldots f(a_{n-1}) = a_n, f(a_n) = a_1. \qquad (35)$$

Indeed: starting from f: $\{a_1,\ldots,a_n\} \to \{b_1,\ldots,b_m\}$, with $a_i > 0$ for i = 1, ..n we have that m ≤ n, otherwise f would not be a function, and for m < n f is not injective, hence $h_f = f \circ f$ is not injective, which is impossible by Lemma 1. Again we have that $\{b_1,\ldots,b_n\} = \{a_1,\ldots,a_n\}$ by $h_f = f \circ f$ and the definition of $h_f$. As in the case of $h_f = f$ we may also here assume that (35) holds: namely that we are in the case of a maximal loop. Indeed, if there were a shorter loop A then we could apply the same reasoning as in the case of $h_f = f$ to $f|_A$ which is possible, by (36), because A is a loop of f, hence also of f∘f and by definition of $h_f$ :

$$h_{f|_A} = (h_f)|_A = (f \circ f)|_A = f|_A \circ f|_A \tag{36}$$

Theorem 10.

$h_f = f \circ f \Leftrightarrow n = 1$ and $f = \{(1,1)\}$

Proof. $h_f = f \circ f$ is equivalent with

$$f(f(f(x))) = x\, f(f(x)) \tag{37}$$

This, combined with (35), yields (for n ≥ 4)

$$a_4 = a_1 a_3,\ a_5 = a_2 a_4,\ \ldots,\ a_n = a_{n-3} a_{n-1},\ a_1 = a_{n-2} a_n,\ a_2 = a_{n-1} a_1,\ a_3 = a_n a_2 \tag{38}$$

This leads to, iteratively, $a_5 = a_2 a_4 = a_1 a_2 a_3$, and finally to

$$a_n = a_1^{E_{n-3}} a_2^{E_{n-4}} a_3^{E_{n-2}} \tag{39}$$

where $(E_n)_n$ is the sequence of numbers obtained in the continuous case. For $a_1, a_2, a_3$ we have:

$$a_1 = a_1^{E_{n-2}} a_2^{E_{n-3}} a_3^{E_{n-1}}$$

$$a_2 = a_1^{E_{n-1}} a_2^{E_{n-2}} a_3^{E_n}$$

$$a_3 = a_1^{E_n} a_2^{E_{n-1}} a_3^{E_{n+1}} \tag{40}$$

It suffices now to prove that the system (40) has only the solution $a_1 = a_2 = a_3 = 1$, and hence $a_i = 1$ for all i = 1, .., n. Then n= 1 and f(1) = 1. Now, setting $x = \ln(a_1)$, $y = \ln(a_2)$ and $z = \ln(a_3)$, the system (40) is equivalent with

$$x = x E_{n-2} + y E_{n-3} + z E_{n-1}$$

$$y = x E_{n-1} + y E_{n-2} + z E_n$$

$$z = x E_n + y E_{n-1} + z E_{n+1} \tag{41}$$

System (41) is a homogeneous system with main determinant:

$$\Delta_n = \begin{vmatrix} E_{n-2} - 1 & E_{n-3} & E_{n-1} \\ E_{n-1} & E_{n-2} - 1 & E_n \\ E_n & E_{n-1} & E_{n+1} - 1 \end{vmatrix} \qquad (42)$$

Now, it suffices to show that $\forall n \in \mathbb{N}: \Delta_n \neq 0$. To prove this we take several steps, each leading to a remarkable result about the sequence $(E_n)_n$.

Step 1. $$\forall n \in \mathbb{Z}: \delta_n =: \begin{vmatrix} E_{n-2} & E_{n-3} & E_{n-1} \\ E_{n-1} & E_{n-2} & E_n \\ E_n & E_{n-1} & E_{n+1} \end{vmatrix} = 1 \qquad (43)$$

Proof. By induction in $\mathbb{Z}$. This means that we will prove $\delta_0 = 1$, $\delta_n = 1 \Rightarrow \delta_{n+1} = 1$, and $\delta_n = 1 \Rightarrow \delta_{n-1} = 1$. Now

$$\delta_0 = \begin{vmatrix} 1 & 0 & 0 \\ 0 & 1 & 0 \\ 0 & 0 & 1 \end{vmatrix} = 1$$

Assume now that $\exists n : \delta_n = 1$ then we consider

$\delta_{n+1} = \begin{vmatrix} E_{n-1} & E_{n-2} & E_n \\ E_n & E_{n-1} & E_{n+1} \\ E_{n+1} & E_n & E_{n+2} \end{vmatrix}$. Now we interchange columns 1 and 2, and next, the new columns 2 and 3. This leads to:

$\delta_{n+1} = \begin{vmatrix} E_{n-2} & E_n & E_{n-1} \\ E_{n-1} & E_{n+1} & E_n \\ E_n & E_{n+2} & E_{n+1} \end{vmatrix}$. Then:

$$\delta_{n+1} - \delta_n = \begin{vmatrix} E_{n-2} & E_n - E_{n-3} & E_{n-1} \\ E_{n-1} & E_{n+1} - E_{n-2} & E_n \\ E_n & E_{n+2} - E_{n-1} & E_{n+1} \end{vmatrix} = \begin{vmatrix} E_{n-2} & E_{n-1} & E_{n-1} \\ E_{n-1} & E_n & E_n \\ E_n & E_{n+1} & E_{n+1} \end{vmatrix} = 0.$$

Hence, if $\delta_n = 1$,

then also $\delta_{n+1} = 1$. For $\delta_{n-1}$ we have:

$$\delta_{n-1} = \begin{vmatrix} E_{n-3} & E_{n-4} & E_{n-2} \\ E_{n-2} & E_{n-3} & E_{n-1} \\ E_{n-1} & E_{n-2} & E_n \end{vmatrix}$$

We interchange rows 1 and 2, and then the new rows 2 and 3. This yields:

$$\delta_{n-1} = \begin{vmatrix} E_{n-2} & E_{n-3} & E_{n-1} \\ E_{n-1} & E_{n-2} & E_n \\ E_{n-3} & E_{n-4} & E_{n-2} \end{vmatrix} \tag{44}$$

From this we obtain:

$$\delta_n - \delta_{n-1} = \begin{vmatrix} E_{n-2} & E_{n-3} & E_{n-1} \\ E_{n-1} & E_{n-2} & E_n \\ E_n - E_{n-3} & E_{n-1} - E_{n-4} & E_{n+1} - E_{n-2} \end{vmatrix} =$$

$\begin{vmatrix} E_{n-2} & E_{n-3} & E_{n-1} \\ E_{n-1} & E_{n-2} & E_n \\ E_{n-1} & E_{n-2} & E_n \end{vmatrix} = 0$. Hence, if δ$_n$ = 1, we also have δ$_{n-1}$ = 1. This proves Step 1 □

Although we proved (43) $\forall n \in \mathbb{Z}$, we will only need the case $n \in \mathbb{N}$ further on.

Step 2. $\forall n \in \mathbb{Z}$:

$$E_n E_{n+3} - E_{n+1} E_{n+2} = E_{-n-4} = E_{n+3}^2 - E_{n+2} E_{n+4} \tag{45}$$

Proof. These equalities are easy to check for the cases n=1,2,3. Then we apply induction from the cases (n-2) and (n-3) to the case for n.

First, we prove the left-hand side equality.

*E$_{n+1}$E$_{n+2}$ + E$_{-n-4}$ = (E$_{n-2}$+E$_n$)(E$_{n-1}$+E$_{n+1}$) + E$_{-n-1}$ - E$_{-n-2}$*

*= E$_{n-2}$E$_{n-1}$+E$_{n-2}$E$_{n+1}$+E$_n$E$_{n-1}$+E$_n$E$_{n+1}$ + E$_{-n-1}$ - E$_{-n-2}$*

Applying the (assumed to be known) equalities *E$_{n-2}$E$_{n-1}$ = E$_{n-3}$E$_n$ − E$_{-n-1}$* and *E$_{n-2}$E$_{n+1}$ = E$_n$E$_{n-1}$ + E$_{-n-2}$* we get:

*= E$_{n-3}$E$_n$ − E$_{-n-1}$ + E$_n$E$_{n-1}$ + E$_{-n-2}$ + E$_n$E$_{n-1}$ + E$_n$E$_{n+1}$+ E$_{-n-1}$ - E$_{-n-2}$*

*= E$_n$(E$_{n-1}$+E$_{n+1}$) + E$_n$(E$_{n-3}$+E$_{n-1}$)*

*=E$_n$E$_{n+2}$ + E$_n$E$_n$ = E$_n$(E$_n$+E$_{n+2}$) = E$_n$ E$_{n+3}$*

This can analogously be proved for $\forall n \in \mathbb{Z}^-$.

Next, we prove the equality on the right-hand side with the same induction hypothesis.

$E_{n+2}E_{n+4} + E_{-n-4} = (E_{n-1}+E_{n+1})(E_{n+1}+E_{n+3}) + E_{-n-1} - E_{-n-2}$

$= E_{n-1}E_{n+1} + E_{n-1}E_{n+3} + (E_{n+1})^2 + E_{n+1}E_{n+3} + E_{-n-1} - E_{-n-2}$

Applying the equalities: $E_{n-1}E_{n+1} = (E_n)^2 - E_{-n-1}$ and $(E_{n+1})^2 - E_{-n-2} = E_n E_{n+2}$

we get:

$= (E_n)^2 + E_{n-1}E_{n+3} + E_n E_{n+2} + E_{n+1}E_{n+3}$

$= E_n(E_n+E_{n+2}) + E_{n+3}(E_{n-1} + E_{n+1})$

$= E_n E_{n+3} + E_{n+3}E_{n+2} = E_{n+3}(E_n+E_{n+2}) = (E_{n+3})^2$

This proves the right-hand side of (45), as the cases n =1,2,3 can be checked directly, and the case $n \in \mathbb{Z}^-$ can be shown similarly as the case $n \in \mathbb{Z}^+$.

Step 3: $\forall n \in \mathbb{Z}$:
$$(E_{n+2})^2 - E_n E_{n+4} = E_{-n-6} = (E_{n+3})^2 - E_n E_{n+6} \qquad (46)$$

Proof of the left-hand side equality:

By (45) we know that

$E_{n+2}E_{n+5} - E_{n+3}E_{n+4} = E_{-n-6}$. Hence we have to show that

$(E_{n+2})^2 - E_n E_{n+4} = E_{n+2}E_{n+5} - E_{n+3}E_{n+4}$ or $E_{n+2}E_{n+5} - E_{n+4}(E_{n+3}-E_n) = (E_{n+2})^2$

This holds as $E_{n+3} = E_n+E_{n+2}$ and $E_{n+5} = E_{n+2} + E_{n+4}$.

Next, we have to prove that:

$(E_{n+2})^2 - E_n E_{n+4} = (E_{n+3})^2 - E_n E_{n+6}$

or: $(E_{n+2}+E_{n+3})(E_{n+3}-E_{n+2}) = E_n(E_{n+6} - E_{n+4})$.

This holds as $E_{n+3}=E_n+E_{n+2}$, hence $(E_{n+2}+E_{n+3})(E_{n+3}-E_{n+2}) = (E_{n+2}+E_{n+3})E_n$ and $E_{n+6}-E_{n+4} = E_{n+2}+ E_{n+3}$, because $E_{n+2} + E_{n+4}= E_{n+5} = E_{n+6}-E_{n+3}$. □

Formulae (45) and (46) make it possible to calculate values of $E_n$, n < 0 from values of $E_n$, n > 0.

Step 4. $\forall n \in \mathbb{Z}$

$$\Delta_n = E_{n+1} - E_{-n+1} + 2(E_{n-2}-E_{-n-2}) \qquad (47)$$

Proof. Determinant $\Delta_n$ (42) is equal to

$$\delta_n - \begin{vmatrix} E_{n-2} - 1 & E_n \\ E_{n-1} & E_{n+1} - 1 \end{vmatrix} - \begin{vmatrix} E_{n-2} & E_{n-1} \\ E_n & E_{n+1} - 1 \end{vmatrix} - \begin{vmatrix} E_{n-2} & E_{n-3} \\ E_{n-1} & E_{n-2} \end{vmatrix}$$

$= \delta_n - ((E_{n-2} -1)(E_{n+1} -1) - E_n E_{n-1} + E_{n-2}(E_{n+1} -1) - E_n E_{n-1} + (E_{n-2})^2 - E_{n-1}E_{n-3}$

$= \delta_n - ( E_{n-2}E_{n+1} - E_{n-2} - E_{n+1} + 1 - E_n E_{n-1} + E_{n-2}E_{n+1} - E_{n-2} - E_n E_{n-1} +(E_{n-2})^2 - E_{n-1}E_{n-3})$

$= \delta_n - ( E_{-n+1} + 2E_{-n-2} - 2E_{n-2} - E_{n+1} + 1)$ ( by (45) with n equal to n-2)

which by (43) proves (47). □

Step 5. $\forall n \in \mathbb{N}$: $E_{n+1} > E_{-n+1}$ and $E_{n-2} > E_{-n-2}$ (48)

Proof. We will use formula (17) for the $(E_n)_n$ :

$$E_n = C_1 a_0^n + C_2 \rho^n \cos(n\zeta) + C_3 \rho^n \sin(n\zeta)$$

with $C_1 = -C_2 \approx 0.4173489$, $C_3 \approx 0.367685$, $a_0 \approx 1.4648493$; $\rho \approx 0.8265432$ and $\zeta \approx 1.8558416$ (radians). Hence:

$$E_{n+1} = C_1 a_0^{n+1} - C_1 \rho^{n+1} \cos((n+1)\zeta) + C_3 \rho^{n+1} \sin((n+1)\zeta)$$

$$E_{-n+1} = \frac{C_1}{a_0^{n-1}} - \frac{C_1}{\rho^{n-1}} \cos((n-1)\zeta) - \frac{C_3}{\rho^{n-1}} \sin((n-1)\zeta)$$

Hence, $E_{n+1} > C_1 a_0^{n+1} - (C_1 + C_3)\rho^{n+1}$ and $E_{-n+1} < \frac{C_1}{a_0^{n-1}} + \frac{C_1+C_3}{\rho^{n-1}}$.

The required inequality $E_{n+1} > E_{-n+1}$ holds if $C_1 a_0^{n+1} - (C_1 + C_3)\rho^{n+1} > \frac{C_1}{a_0^{n-1}} + \frac{C_1+C_3}{\rho^{n-1}}$ or equivalently, if:

$$C_1 \left(a_0^{n+1} - \frac{1}{a_0^{n-1}}\right) > (C_1 + C_3)\left(\rho^{n+1} + \frac{1}{\rho^{n-1}}\right)$$

This inequality holds if

$$C_1 a_0^{n+1} - (C_1 + C_3)\frac{1}{\rho^{n-1}} > 2C_1 + C_3$$

or

$$(0.4173489)(1.4648493)^{n+1} - (0.7850339)(1.2098581)^{n-1} > 1.2023829$$

$$(49).$$

This inequality holds from all $n \geq n_0$. Direct calculation shows that $n_0 = 4$. Moreover, we can readily show that $E_{n+1} > E_{-n+1}$ for $n = 1,2,3$ (see Table 1) so that the left-hand part of (48) is proved for all $n \in \mathbb{N}$. Now we prove the right-hand side of (48).

Similarly (as before) we have that

$$E_{n-2} > C_1 a_0^{n-2} - (C_1 + C_3)\rho^{n-2}$$

and $E_{-n-2} < \frac{C_1}{a_0^{n+2}} + (C_1 + C_3)\frac{1}{\rho^{n+2}}$.

So, $E_{n-2} > E_{-n-2}$ holds if $C_1 \left(a_0^{n-2} - \frac{1}{a_0^{n+2}}\right) > (C_1 + C_3)\left(\rho^{n-2} + \frac{1}{\rho^{n+2}}\right)$.

For this, it suffices to prove that

$$C_1 a_0^{n-2} - (C_1 + C_3)\frac{1}{\rho^{n+2}} > 2C_1 + C_3$$

or:

$(0.4173489)(1.4648493)^{n-2} - (0.7850339)(1.2098581)^{n+2} > 1.2023829$

(50)

Inequality (50) holds for all n > 10. Moreover, one can readily check (see Table 1) that $E_{n-2} > E_{-n-2}$ holds for n = 1, 2, …, 10, proving (48).

Step 6. This concludes the proof that $\forall n \in \mathbb{N}, n \geq 4: \Delta_n > 0,$ hence $\neq 0$. Finally, for n =1,2,3, we argue directly on f. For n = 1, we trivially find f = {(1,1)}, using (37) (since $a_1 \neq 0$). For n =2 we have, by (37) that $a_2 = a_1 a_2$ and $a_1 = a_2^2$, whence $a_1 = a_2 = 1$ and f = {(1,1)}. For n = 3 , (37) gives $a_1 = a_1 a_3$ , $a_2 = a_2 a_1$ and $a_3 = a_3 a_2$, hence (since $a_i \neq 0$), $a_1 = a_2 = a_3 = 1$ and f = {(1,1)}. □

## 7. The problem $h_{g \circ g} = g$, with f = g∘g in the discrete case

As in previous cases, we may assume, now for the case g, that

$g(a_1) = a_2, g(a_2) = a_3, …, g(a_{n-1}) = a_n , g(a_n) = a_1.$ (51)

Indeed, given g: {$a_1,…,a_n$} → {$b_1, …, b_m$} it follows that m = n, because of the injectivity of g and hence also {$a_1,…,a_n$} = {$b_1, …, b_n$} by the definition of $h_{g \circ g} = g$. Similar to the case $h_f = f \circ f$ we may assume that we have a maximal loop in (51). If there would exist a shorter loop A then we apply the earlier reasoning on $g|_A$, which is allowed because $h_{g|_A \circ g|_A} = h_{(g \circ g)|_A} = h_{g \circ g}|_A = g|_A$, as A is a loop of g and by definition of h. If g = {(1,1)} then f = g∘g = {(1,1)}.

Theorem 11. $h_{g \circ g} = g$, with f = g∘g ⇔ n =1 and f={(1,1)}

30Proof. The proof is analogous to that in the case $h_f = f \circ f$. Hence we will only provide a sketch, including the corresponding formulae. We will show that $a_1 = \ldots = a_n = 1$ and hence $g = \{(1,1)\}$.

The equality $h_{g \circ g} = g$ is equivalent with

$$g(g(g(x))) = x\, g(x) \tag{52}$$

From this, we derive (using the sequence $(D_n)_n$ introduced in the continuous case)

$$a_n = a_1^{D_{n-1}} a_2^{D_{n+1}} a_3^{D_n} \quad (n \geq 4) \tag{53}$$

and

$$a_1 = a_1^{D_n} a_2^{D_{n+2}} a_3^{D_{n+1}}$$
$$a_2 = a_1^{D_{n+1}} a_2^{D_{n+3}} a_3^{D_{n+2}}$$
$$a_3 = a_1^{D_{n+2}} a_2^{D_{n+4}} a_3^{D_{n+3}} \tag{54}$$

With $x = \ln(a_1)$, $y = \ln(a_2)$ and $z = \ln(a_3)$ system (54) leads to:

$$x = x\, D_n + y\, D_{n+2} + z\, D_{n+1}$$
$$y = x\, D_{n+1} + y\, D_{n+3} + z\, D_{n+2}$$
$$z = x\, D_{n+2} + y\, D_{n+4} + z\, D_{n+3} \tag{55}$$

The system (55) is homogeneous. It has only the 0-solution, hence $a_1 = a_2 = a_3 = 1$, resulting in $a_i = 1$ for all $i = 1, \ldots, n$ iff

$$\Delta'_n = \begin{vmatrix} D_n - 1 & D_{n+2} & D_{n+1} \\ D_{n+1} & D_{n+3} - 1 & D_{n+2} \\ D_{n+2} & D_{n+4} & D_{n+3} - 1 \end{vmatrix} \neq 0 \tag{56}$$

We use the notation $\delta'_n$ for the determinant

$$\begin{vmatrix} D_n & D_{n+2} & D_{n+1} \\ D_{n+1} & D_{n+3} & D_{n+2} \\ D_{n+2} & D_{n+4} & D_{n+3} \end{vmatrix} \tag{57}$$

Now we prove:

(I) $$\forall n \in \mathbb{Z}: \delta'_n = 1 \tag{58}$$

Proof. We will not include all details, but we will show that the reasoning to prove (58) is the same as in the ($E_n$)-case. We first note that clearly $\delta'_1 = 1$ and will show that if $\delta'_n = 1$ it follows that $\delta'_{n+1} = 1$ ( $n \in \mathbb{N}$).

$$\delta'_{n+1} = \begin{vmatrix} D_{n+1} & D_{n+3} & D_{n+2} \\ D_{n+2} & D_{n+4} & D_{n+3} \\ D_{n+3} & D_{n+5} & D_{n+4} \end{vmatrix} = \begin{vmatrix} D_{n+3} & D_{n+2} & D_{n+1} \\ D_{n+4} & D_{n+3} & D_{n+2} \\ D_{n+5} & D_{n+4} & D_{n+3} \end{vmatrix}$$

From this, it follows that

$$\delta'_{n+1} - \delta'_n = \begin{vmatrix} D_{n+3} - D_n & D_{n+2} & D_{n+1} \\ D_{n+4} - D_{n+1} & D_{n+3} & D_{n+2} \\ D_{n+5} - D_{n+2} & D_{n+4} & D_{n+3} \end{vmatrix} = \begin{vmatrix} D_{n+1} & D_{n+2} & D_{n+1} \\ D_{n+2} & D_{n+3} & D_{n+2} \\ D_{n+3} & D_{n+4} & D_{n+3} \end{vmatrix} = 0$$

(where we have used equation (33)). A similar inductive reason, from n to (n-1) proves the result $\forall n \in \mathbb{Z}$. □

Step (II):

$$\forall n \in \mathbb{Z}: D_n^2 - D_{n-1}D_{n+1} = D_{-n+3} = D_n D_{n+3} - D_{n+1}D_{n+2} \quad (59)$$

Step (III): $\forall n \in \mathbb{Z}: \Delta'_n = D_n - D_{-n} + 2(D_{n+3} - D_{-n+3})$ (60)

Step (IV): $\forall n \in \mathbb{N}$

$$D_n > D_{-n} \text{ and } D_{n+3} > D_{-n+3} \quad (61)$$

From (61) it follows that $\Delta'_n \neq 0$ (actually > 0)

Proof of (61). We have that $\forall n \in \mathbb{Z}$, by (34):

$D_n = C'_1 c_1^n + C'_2 (\rho')^n \cos(n\zeta') + C'_3 (\rho')^n \sin(n\zeta')$, with $D_0 = 1$, $D_1 = D_2 = 0$.

We know already that

$c = c_1 \approx 1.3247178$

$c_{2,3} \approx -0.6623589 \pm i\, 0.5622796$

from which we derive that $\rho' \approx 0.868837$ and $\zeta \approx 2.4377348$ rad

Finally, we have: $C_1' \approx 0.1770089$, $C_2' = 1 - C_1' \approx 0.8229911$, $C_3' \approx 0.5524453$. From this, we obtain (61) in a similar way as in the case for $(E_n)_n$. □

For further use, we present the following table of $D_n$ values.

Table 3. D-numbers

| n   | -1  | 0   | 1   | 2   | 3   | 4   | 5   | 6   | 7   | 8   | 9   | 10  | 11  | 12  | 13  |
|-----|-----|-----|-----|-----|-----|-----|-----|-----|-----|-----|-----|-----|-----|-----|-----|
| $D_n$ | -1  | 1   | 0   | 0   | 1   | 0   | 1   | 1   | 1   | 2   | 2   | 3   | 4   | 5   | 7   |
| n   | -16 | -15 | -14 | -13 | -12 | -11 | -10 | -9  | -8  | -7  | -6  | -5  | -4  | -3  | -2  |
| $D_n$ | -3  | 7   | -7  | 4   | 0   | -3  | 4   | -3  | 1   | 1   | -2  | 2   | -1  | 0   | 1   |

## 8. Relations between the D, E, and F-sequences

As $D_n$, $E_n$, and $F_n \geq 0$ for $n = 0,1,2$, it follows by definition that $\forall n \in \mathbb{N}$, $D_n$, $E_n$, and $F_n \geq 0$. These inequalities do not hold for $n \in \mathbb{Z}^-$. Again, by their very definitions, these sequences are strictly increasing. Differences, for the same value of n, between these sequences, increase. This follows from the following results.

Theorem 12. For $n \in \mathbb{N}$,

$$D_n < E_{n-1} \ (n \geq 4) \tag{62}$$

$$E_n < F_{n-1} \ (n \geq 6) \tag{63}$$

Proof. By induction. For (62) we assume that the inequality holds for n-2 and n-3 and then provide proof for n. By assumption, we know that

$$D_n = D_{n-3} + D_{n-2} < E_{n-4} + E_{n-3} \leq E_{n-4} + E_{n-2} = E_{n-1}$$

One can then check numerically (see tables) that (62) holds for $n \geq 4$.

For (63) we assume that the inequality holds for n-1 and n-3 and then prove for n. We then know that $E_n = E_{n-3} + E_{n-1} < F_{n-4} + F_{n-2} \leq F_{n-3} + F_{n-2} = F_{n-1}$. Again one has to check numerically that (63) holds for n ≥ 6. □

Corollary. For $n \in \mathbb{N}$ we have

$$D_{n+1} - D_n < E_{n+1} - E_n, \quad n \geq 5 \qquad (64)$$

$$E_{n+1} - E_n < F_{n+1} - F_n, \quad n \geq 4 \qquad (65)$$

Proof. By Theorem 12 we know that

$D_{n+1} - D_n \leq D_{n+2} - D_n = D_{n-1} < E_{n-2} = E_{n+1} - E_n$ for n ≥ 5. Further, by this same theorem, we know that

$E_{n+1} - E_n = E_{n-2} < F_{n-3} < F_{n-1} = F_{n+1} - F_n$ for n ≥ 8, we can easily check that this inequality also holds for n = 4,5,6 and 7. □

Corollary. The sequences $(E_n - D_n)_n$ and $(F_n - E_n)_n$ are strictly increasing for n ≥ 5 and tend to ∞ for n tending to ∞:

$$\lim_{n \to +\infty}(E_n - D_n) = \lim_{n \to +\infty}(F_n - E_n) = +\infty \qquad (66)$$

Proof. This result follows immediately from (64) and (65), observing that, for $n \in \mathbb{N}$, (64) and (65) imply that

$D_{n+1} - D_n + 1 \leq E_{n+1} - E_n$ and $E_{n+1} - E_n + 1 \leq F_{n+1} - F_n$ □

More general than (66) we have the following theorem.

Theorem 13. $\forall k \in \mathbb{N} \cup \{0\}$, k fixed we have:

$$\lim_{n \to +\infty}(E_{n-k} - D_n) = \lim_{n \to +\infty}(F_{n-k} - E_n) = +\infty \qquad (67)$$

Proof. $\forall k \in \mathbb{N} \cup \{0\}$, we apply the goniometric form of $E_{n-k}$ and of $D_n$, see (17) and (34). Then $E_{n-k} - D_n = C_1 a_0^{n-k} + C_2 \rho^{n-k} \cos((n-k)\zeta) + C_3 \rho^{n-k} \sin((n-k)\zeta) - (C_1' c_1^n + C_2'(\rho')^n \cos(n\zeta') + C_3'(\rho')^n \sin(n\zeta'))$

$\geq C_1 a_0^{n-k} - (C_2 + C_3)\rho^{n-k} - C_1' c_1^n - (C_2' + C_3')(\rho')^n =$

$$\frac{C_1}{a_0^k} a_0^n - C_1' c_1^n - [(C_2 + C_3)\rho^{n-k} + (C_2' + C_3')(\rho')^n].$$

Now, $0 < \rho, \rho' < 1$, and hence the term between square brackets tends to zero for $n \to +\infty$. The first two terms tend to $+\infty$ as $1 < C_1 < a_0$ and (with k fixed):

$$\frac{C_1}{a_0^k} a_0^n - C_1' C_1^n = C_1^n \left( \frac{C_1}{a_0^k} \left(\frac{a_0}{C_1}\right)^n - C_1' \right)$$

We further have: $F_{n-k} - E_n$

$$= \frac{\alpha^{n-k} - \beta^{n-k}}{\sqrt{5}} - (C_1 a_0^n + C_2 \rho^n \cos(n\zeta) + C_3 \rho^n \sin(n\zeta))$$

with $\alpha = \frac{\sqrt{5}+1}{2} > 1$, and $-1 < \beta = \frac{1-\sqrt{5}}{2} < 0$.

Hence, with k fixed, $\lim_{n \to +\infty}(F_{n-k} - E_n) = \frac{1}{\sqrt{5}} \lim_{n \to +\infty} \alpha^{n-k} = +\infty$ □

We note that we provided two different proofs for the case k = 0.

Corollary. $\forall k \in \mathbb{N} \cup \{0\}$

(i) $\exists n_0(k) \in \mathbb{N}$ such that for $n \geq n_0(k)$

$$D_n < E_{n-k} \tag{68}$$

(ii) $\exists n_1(k) \in \mathbb{N}$ such that for $n \geq n_1(k)$

$$E_n < F_{n-k} \tag{69}$$

For the reader's information, we include a table of minimal values for $n_0(k)$ and $n_1(k)$, noting that it includes (62) and (63).

Table 4. A table of minimal values for $n_0(k)$ and $n_1(k)$

| k | n₀(k) | n₁(k) |
|---|---|---|
| 0 | 4 | 3 |
| 1 | 4 | 6 |
| 2 | 6 | 10 |
| 3 | 7 | 14 |
| 4 | 9 | 19 |

Corollary

For $\forall k \in \mathbb{N} \cup \{0\}$ fixed, the following inequalities hold for sufficiently large values of n:

$$D_{n+k-1} - D_{n+k} < E_{n+1} - E_n \qquad (70)$$

$$E_{n+k-1} - E_{n+k} < F_{n+1} - F_n \qquad (71)$$

Proof. $D_{n+k-1} - D_{n+k} \leq D_{n+k+2} - D_{n+k} = D_{n+k-1} < E_{n-2} = E_{n+1} - E_n$, by using (68) for k+1. We can prove (71) in a similar way, using (69). □

Corollary.

$$\lim_{n \to +\infty}(F_n - (E_n + D_n)) = +\infty \qquad (72)$$

Proof. For $n \in \mathbb{N}$, large, we have: $F_n - (D_n + E_n) = F_{n-1} + F_{n-2} - (D_{n-2} + D_{n-3} + E_{n-1} + E_{n-3}) = F_{n-3} + F_{n-2} + F_{n-4} + F_{n-3} - D_{n-2} - D_{n-3} - E_{n-1} - E_{n-3}$.

By (71) we know that $F_{n-4} - E_{n-3} > F_{n-5} - E_{n-4}$ and hence

$$F_{n-4} - E_{n-3} \geq 1 + F_{n-5} - E_{n-4} \qquad (73)$$

Applying this reasoning three times (using (70) and (71)) yields:

$F_n - (D_n + E_n) \geq 4 + F_{n-5} - F_{n-4} + F_{n-4} - E_{n-2} + F_{n-3} - D_{n-3} + F_{n-4} - D_{n-4}$.

Applying again (70) and (71) yields (72). □

Also here we have that the sequence $(F_n - (E_n + D_n))_n$ is strictly increasing in n. This follows from the fact that we can replace $E_n$ in (69) by $(E_n + D_n)$ and similarly for (71) where we replace ($E_{n+1} - E_n$) by ($E_{n+1}+D_{n+1}$)-($E_n+D_n$) (a similar proof can be given), from which follows the strict increase of ($F_n - (E_n+D_n)$)$_n$. □

Finally, we mention some further properties. If we define $\forall\, n \in \mathbb{Z}$: $G_n = E_n + F_n$ and $H_n = D_n+E_n+F_n$, then $G_n = H_n - F_n$. It is now easy to show that $\forall n \in \mathbb{N}$:

$$G_n + G_{n+1} \leq G_{n+3} \leq G_n + G_{n+2} \tag{74}$$

$$H_n + H_{n+1} \leq H_{n+3} \not\leq H_n + H_{n+2} \tag{75}$$

and strictly smaller from a small n ∈ ℕ on, the case $\not\leq$ happens e.g., for n = 2,4, …

## 9. Further study of the characteristic determinant of the sequences $(F_n)_n$, $(E_n)_n$, and $(D_n)_n$

Concretely, we will study and generalize formulae (6), (43), and (57). As a motivation for further developments, we begin with the $(F_n)_n$ sequence, followed by a thorough study of the $(E_n)_n$ case and an outline of how the $(D_n)_n$ case can be handled.

9.1 The $(F_n)_n$ case (Fibonacci)

We know that $\forall\, n \in \mathbb{Z}: F_{n-1}F_{n+1} - F_n^2 = (-1)^n$ (6)

We define for m, n ∈ ℤ

$$\Delta_{n,m} = \begin{vmatrix} F_n & F_{n+1} \\ F_m & F_{m+1} \end{vmatrix} \tag{76}$$

This leads to the following theorem, known as Vajda's identity.

Theorem 14. Vajda's identity (Vajda, 2008)

$$\forall\, m, n \in \mathbb{Z}: \Delta_{n,m} = F_{m-n}(-1)^{n+1} = F_{n-m}(-1)^m \tag{77}$$

Proof. As $\Delta_{m,m} = 0$ and $\Delta_{n,m} = -\Delta_{m,n}$ it suffices to use an induction argument in one variable. We will show that if (77) holds for n-1 and n $\in \mathbb{N}$ then (77) holds for n+1.

$\Delta_{n+1,m}$ = $F_{n+1}F_{m+1} - F_m F_{n+2}$ = $(F_n + F_{n-1})F_{m+1} - F_m(F_n + F_{n+1})$

= $F_n F_{m+1} - F_m F_{n+1} + F_{n-1} F_{m+1} - F_m F_n$ = $F_{m-n}(-1)^{n+1} + F_{m-n+1}(-1)^n$

= $(F_{m-n} - F_{m-n+1})(-1)^{n+1}$ = $- F_{m-n-1}(-1)^{n+1}$ = $F_{m-n-1}(-1)^{n+2}$

This is the left-hand side of (77). We further have by (3):

$$F_{m-n} = F_{n-m}(-1)^{(n-m+1)}$$

Hence $\Delta_{n,m} = F_{m-n}(-1)^{n+1}$ implies that $\Delta_{n,m} = F_{n-m}(-1)^{n-m+1+n+1}$

= $F_{n-m}(-1)^{-m}$ = $F_{n-m}(-1)^m$ □

We observe that formula (77) generalizes formula (8) (take m = n+1). It is moreover the key to the solution to the following problem.

Problem. Given a Fibonacci sequence, i.e., a sequence such that each term is the sum of the two previous ones, (F'$_n$)$_n$ with F'$_1$ = a and F'$_2$ = b (a ≠ b), determine a and b if two values (and their different indices) are given.

We first prove a lemma.

Lemma 3.
$$\forall n \in \mathbb{Z}: F'_{n+2} = a\, F_n + b\, F_{n+1} \tag{78}$$

Proof. This is easy to see as $F'_1 = a, F'_2 = b, F'_3 = a+b, F'_4 = 2a+b, \ldots$ Formally, using induction assuming that the cases n-1 and n are known, we prove the formula for n+1.

$$F'_{n+3} = F'_{n+1} + F'_{n+2} = a\, F_{n-1} + b\, F_n + a F_n + b\, F_{n+1}$$
$$= a(F_{n-1} + F_n) + b(F_n + F_{n+1}) = a\, F_{n+1} + b F_{n+2}\ .$$

Similarly, assuming known the case n+1 and n+2 we prove the case n-1:

$$F'_{n+1} = F'_{n+3} - F'_{n+2} = a\,F_{n+1} + b\,F_{n+2} - aF_n - b\,F_{n+1}$$
$$= a(F_{n+1} - F_n) + b(F_{n+2} - F_{n+1}) = a\,F_{n-1} + bF_n \quad \square$$

Denoting the two given values in the $(F'_n)_n$ sequence by $F'_{n+2}$ and $F'_{m+2}$ (m≠n) then we have, by (78),

$$F'_{n+2} = a\,F_n + b\,F_{n+1}$$

$$F'_{m+2} = a\,F_m + b\,F_{m+1}$$

From which we find:

$$a = \frac{\begin{vmatrix} F'_{n+2} & F_{n+1} \\ F'_{m+2} & F_{m+1} \end{vmatrix}}{\begin{vmatrix} F_n & F_{n+1} \\ F_m & F_{m+1} \end{vmatrix}}, \quad b = \frac{\begin{vmatrix} F_n & F'_{n+2} \\ F_m & F'_{m+2} \end{vmatrix}}{\begin{vmatrix} F_n & F_{n+1} \\ F_m & F_{m+1} \end{vmatrix}},$$

With (45) this becomes:

$$a = \frac{(-1)^m}{F_{n-m}} (F'_{n+2} F_{m+1} - F'_{m+2} F_{n+1})$$

$$b = \frac{(-1)^m}{F_{n-m}} (F'_{m+2} F_n - F'_{n+2} F_m)$$

This answers the problem (note that $F_{n-m} \neq 0$ as n ≠ m).

An example using the Lucas sequence: 1,3,4,7,11,18,29,47,76, …

Assume $F'_6 = 18$ and $F'_7 = 29$ are given, i.e., n+2 = 6 and m+2 + 7. Then:

a = $(-1)^5$ ((18).(8) − (29).(5)) = 1

b = $(-1)^5$ ((29).(3) − (18).(5)) = 3.

The given F'- values do not have to be consecutive ones. For example, if it is given that F'$_6$ = 18 and F'$_8$ = 47, then n+2 = 6 and m+2 = 8 and hence a = (-1)$^6$(18x13 − 47 x 5) = 1 and b = (-1)$^6$ (47x3 − 18x8) = 3

We next note the following, simple result.

Proposition 4

The Fibonacci sequence (F'$_n$)$_n$ is the right-hand side of the classical Fibonacci sequence (F$_n$)$_n$ if and only if $\exists\, n_0 \in \mathbb{Z}$ such that $a =: F'_1 = F_{n_0}$ and $b =: F'_2 = F_{n_0+1}$

Proof. The implication from left to right is trivial

Assuming now that $\exists\, n_0 \in \mathbb{Z}$ such that $a =: F'_1 = F_{n_0}$ and $b =: F'_2 = F_{n_0+1}$ then it follows from the basic Fibonacci property $F'_n + F'_{n+1} = F'_{n+2}$ that (F'$_n$)$_n$ is the right-hand side of the classical Fibonacci sequence.

9.2 The sequence (E$_n$)$_n$

We already know from (43) that

$$\forall\, n \in \mathbb{Z}: \delta_n =: \begin{vmatrix} E_{n-2} & E_{n-3} & E_{n-1} \\ E_{n-1} & E_{n-2} & E_n \\ E_n & E_{n-1} & E_{n+1} \end{vmatrix} = 1$$

We will now extend the determinant used in (43) to:

$$\Delta_{n,m,k} =: \begin{vmatrix} E_n & E_{n-1} & E_{n+1} \\ E_m & E_{m-1} & E_{m+1} \\ E_k & E_{k-1} & E_{k+1} \end{vmatrix} \qquad (79)$$

We see that with m=n+1 and k=m+1, $\Delta_{n,m,k} = \delta_{n+2} = 1$. We can rewrite the variables n, m, and k as n, m = n+i and k = m+j = n+i+j, leading to (already assuming that this determinant is independent of n, see next theorem):

$$\Delta_{i,j} := \Delta_{n,m,k} = \begin{vmatrix} E_n & E_{n-1} & E_{n+1} \\ E_{n+i} & E_{n+i-1} & E_{n+i+1} \\ E_{n+i+j} & E_{n+i+j-1} & E_{n+i+j+1} \end{vmatrix} \tag{80}$$

Theorem 15. $\Delta_{n,n+i,n+i+j}$ is independent of $n \in \mathbb{Z}$ and

$$\Delta_{i,j} = E_i E_{i+j-1} - E_{i+j} E_{i-1} \tag{81}$$

Hence $\Delta_{i,j}$ only depends on m-n = i and k-n = i+j.

Proof. From (80) for n+1 we have:

$$\Delta_{n+1,n+i+1,n+i+j+1} = \begin{vmatrix} E_{n+1} & E_n & E_{n+2} \\ E_{n+i+1} & E_{n+i} & E_{n+i+2} \\ E_{n+i+j+1} & E_{n+i+j} & E_{n+i+j+2} \end{vmatrix}$$

We apply the following equalities:

$E_{n+2} = E_{n-1} + E_{n+1}$

$E_{n+i+2} = E_{n+i-1} + E_{n+i+1}$

$E_{n+i+j+2} = E_{n+i+j-1} + E_{n+i+j+1}$

leading to

$$\Delta_{n+1,n+i+1,n+i+j+1} = \begin{vmatrix} E_{n+1} & E_n & E_{n-1} \\ E_{n+i+1} & E_{n+i} & E_{n+i-1} \\ E_{n+i+j+1} & E_{n+i+j} & E_{n+i+j-1} \end{vmatrix} = \Delta_{n,n+i,n+i+j}$$

after an exchange of the first and the third column, followed by an exchange of the new columns one and two. From this, it follows that in (80) we may take n=0. Hence:

$$\Delta_{i,j} := \Delta_{n,n+i,n+i+j} = \Delta_{n,m,k} = \begin{vmatrix} 0 & 0 & 1 \\ E_i & E_{i-1} & E_{i+1} \\ E_{i+j} & E_{i+j-1} & E_{i+j+1} \end{vmatrix}$$

from which (81) follows. □

From the previous result, we can represent $\Delta$ in two different ways as a matrix, see Tables 5 and 6.

Table 5. $\Delta_{i,j} := \Delta_{n,n+i,n+i+j}$ with $i,j \in \mathbb{N}$.

| i\j | 1 | 2 | 3 | 4 | 5 | 6 | 7 | 8 | … |
|---|---|---|---|---|---|---|---|---|---|
| 1 | 1 | 1 | 1 | 2 | 3 | 4 | 6 | 9 | |
| 2 | 0 | -1 | -1 | -1 | -2 | -3 | -4 | -6 | |
| 3 | -1 | -1 | -1 | -2 | -3 | -4 | -6 | -9 | |
| 4 | 1 | 2 | 2 | 3 | 5 | 7 | 10 | 15 | |
| 5 | 1 | 0 | 0 | 1 | 1 | 1 | 2 | 3 | |
| 6 | -2 | -3 | -3 | -5 | -8 | -11 | -16 | -24 | |
| 7 | 0 | 2 | 2 | 2 | 4 | 6 | 8 | 12 | |
| 8 | 3 | 3 | 3 | 6 | 9 | 12 | 18 | 27 | |
| … | | | | | | | | | |

Table 6. $\Delta_{n,m,k}$ with $m, m \in \mathbb{N}$

| m\k | 1 | 2 | 3 | 4 | 5 | 6 | 7 | 8 | 9 | 10 | 11 | 12 | 13 | … |
|---|---|---|---|---|---|---|---|---|---|---|---|---|---|---|
| 1 | 0 | 1 | 1 | 1 | 2 | 3 | 4 | 6 | 9 | … | | | | |
| 2 | | 0 | 0 | -1 | -1 | -1 | -2 | -3 | -4 | -6 | … | | | |
| 3 | | | 0 | -1 | -1 | -1 | -2 | -3 | -4 | -6 | -9 | … | | |
| 4 | | | | 0 | 1 | 2 | 2 | 3 | 5 | 7 | 10 | 15 | … | |
| 5 | | | | | 0 | 1 | 0 | 0 | 1 | 1 | 1 | 2 | 3 | … |
| 6 | | | | | | 0 | -2 | -3 | -3 | -5 | -8 | -11 | -16 | … |
| 7 | | | | | | | 0 | 0 | 2 | 2 | 2 | 4 | 6 | … |
| 8 | | | | | | | | 0 | 3 | 3 | 3 | 6 | 9 | … |
| … | | | | | | | | | | | | | | |

Further, we have, in Table 6, for m > k: $\Delta_{n,m,k} = -\Delta_{n,k,m}$.

Tables 5 and 6 are both n-independent. The $j^{th}$ column of Table 5 is the $j^{th}$ oblique line (slope equal to -1) in Table 6. We stress that e.g., the value 1 in Table 5 corresponding to the coordinates (1,1), can be found in Table 6 with coordinates (1,2) and is formula (43), there denoted as $\delta_n = 1$.

9.3 Properties of $\Delta_{i,j}$

9.31. We begin this section with Theorm 16.

Theorem 16.

(a) For all $i \in \mathbb{Z}$, fixed, $(\Delta_{i,j})_{j \in \mathbb{N}}$ is an E'-sequence in the sense that it satisfies the relation

$$\Delta_{i,j-3} + \Delta_{i,j-1} = \Delta_{i,j} \qquad (82)$$

with initial values

$$\begin{aligned}
a' &= \Delta_{i,1} = E_i^2 - E_{i-1}E_{i+1} \\
b' &= \Delta_{i,2} = E_i E_{i+1} - E_{i-1}E_{i+2} \\
c' &= \Delta_{i,3} = E_i E_{i+2} - E_{i-1}E_{i+3}
\end{aligned} \qquad (83)$$

where b'= c'

(b) For all j ∈ ℤ, fixed, $(\Delta_{i,j})_{i \in \mathbb{N}}$ is a reverse E'-sequence in the sense that is satisfies the relation

$$\Delta_{i+2,j} + \Delta_{i,j} = \Delta_{i-1,j} \qquad (84)$$

(hence, $E'_{-i-2} + E'_{-i} = E'_{-i+1}$ with $E'_i = \Delta_{-i,j}$) with reverse initial values:

$$a'' = E_j \; ; \; b'' = -E_{j-1} \; ; \; c'' = -E_j \qquad (85)$$

Proof.

(a) From (81) we have:

$$\Delta_{i,j-3} + \Delta_{i,j-1} = E_i E_{i+(j-3)-1} - E_{i+(j-3)}E_{i-1} + E_i E_{i+(j-1)-1} - E_{i+(j-1)}E_{i-1}$$

$$= E_i(E_{i+j-4} + E_{i+j-2}) - E_{i-1}(E_{i+j-3} + E_{i+j-1})$$

$$= E_i E_{i+j-1} - E_{i-1}E_{i+j} = \Delta_{i,j}$$

The equations (83) follow from (81). Further we see that b'= c' because

$$c' = E_i E_{i+2} - E_{i-1}E_{i+3} = E_i(E_{i-1} + E_{i+1}) - E_{i-1}(E_i + E_{i+2})$$

$$= E_i E_{i+1} - E_{i-1}E_{i+2} = b'.$$

(b) $\Delta_{i-1,j} - \Delta_{i,j} = E_{i-1}E_{i+j-2} - E_{i+j-1}E_{i-2} - E_i E_{i+j-1} + E_{i+j}E_{i-1}$

$$= E_{i-1}(E_{i+j-2} + E_{i+j}) - E_{i+j-1}(E_{i-2} + E_i)$$

$$= E_{i-1}E_{i+j+1} - E_{i+j-1}E_{i+1}$$

$\Delta_{i+2,j} = E_{i+2}E_{i+j+1} - E_{i+j-2}E_{i+1}$

$= E_{i+2}E_{i+j+1} - (E_{i+j-1} + E_{i+j+1})E_{i+1}$

$= (E_{i+2} - E_{i+1})E_{i+j+1} - E_{i+j-1}E_{i+1} = E_{i-1}E_{i+j+1} - E_{i+j-1}E_{i+1} = \Delta_{i-1,j} - \Delta_{i,j}$ □

This proves (84). Finally, (85) follows for i=1,2,3 in (81), using the values of $E_0$, $E_1$, $E_2$ and $E_3$:

$$a' = E_1 E_j - E_{j+1} E_0 = E_j$$
$$b' = E_2 E_{j+1} - E_{j+2} E_1 = E_{j+1} - E_{j+2} = -E_{j-1}$$
$$c' = E_3 E_{j+2} - E_{j+3} E_2 = E_{j+2} - E_{j+3} = -E_j$$

The reader may compare these results with Table 5.

### 9.3.2. Some remarkable identities, proved using (81)

(a). $\forall i \in \mathbb{Z}: \Delta_{i,2} = \Delta_{i,3}$  (86)

Proof. $\Delta_{i,3} = E_i E_{i+2} - E_{i+3} E_{i-1} = E_i E_{i+2} - (E_i + E_{i+2})E_{i-1}$
$= E_i(E_{i+2} - E_{i-1}) - E_{i+2} E_{i-1} = E_i E_{i+1} - E_{i+2} E_{i-1} = \Delta_{i,2}$ □

(b$_1$) $\forall i \in \mathbb{Z}: \Delta_{i,1} = E_{-i-1}$  (87)

(b$_2$). $\forall i \in \mathbb{Z}: \Delta_{i,2} (= \Delta_{i,3}) - E_{-i-3}$  (88)

Proof. $\Delta_{i,1} = E_i^2 - E_{i-1} E_{i+1} = E_{-(i-3)-4} = E_{-i-1}$, using the right-hand side of (45).

$\Delta_{i,2} = E_i E_{i+1} - E_{i+2} E_{i-1} = -E_{-(i-1)-4} = -E_{-i-3}$ using the left-hand side of (45). □

Consequently, by (88), (87), and (86) we have:

$$\Delta_{i,1} = -\Delta_{i-2,2} = -\Delta_{i-2,3} \tag{89}$$

(c) $\forall j \in \mathbb{Z}: \Delta_{2,j} = \Delta_{3,j-1}$  (90)

Proof. On the one hand, we have: $\Delta_{3,j-1} = E_3 E_{j+1} - E_{j+2} E_2 = E_{j+1} - E_{j+2}$

while on the other, $\Delta_{2,j} = E_2 E_{j+1} - E_{j+2} E_1 = \Delta_{3,j-1}$ □

(d). $\forall j \in \mathbb{Z}: \Delta_{1,j} = -\Delta_{3,j}$  (91)

Proof. On the one hand, we have: $\Delta_{1,j} = E_1 E_j - E_{j+1} E_0 = E_j$

While on the other: $-\Delta_{3,j} = -(E_3 E_{j+2} - E_{j+3} E_2) = -E_{j+2} + E_{j+3} = E_j$ □

Consequently, by (90) and (91) we have:

$$\Delta_{2,j} = -\Delta_{1,j-1} \tag{92}$$

(e) $\forall j \in \mathbb{Z}: \Delta_{1,j} = \Delta_{5,j+3} = E_j$ (93)

Proof. $\Delta_{5,j+3} = E_5 E_{j+7} - E_{j+8} E_4 = 3 E_{j+7} - 2(E_{j+5} + E_{j+7}) = E_{j+7} - 2E_{j+5}$
$= E_{j+4} + E_{j+6} - E_{j+5} - E_{j+5} = E_{j+4} + E_{j+3} - E_{j+5} = -E_{j+2} + E_{j+3} = E_j = \Delta_{1,j}$ □

In the next step, we will consider, using $\Delta_{i,j}$, the analogous problem for E-sequences, as we solved for the F (Fibonacci) sequences.

Problem. Given an E' sequence $(E'_n)_n$ with initial values $a' = E'_1$, $b' = E'_2$, and $c' = E'_3$, determine *a'*, *b'*, and *c'* when three values from the $(E'_n)_n$-sequence are given (more details about these three values follow further). We begin with a lemma.

Lemma 4. $\forall n \in \mathbb{Z}$:

$$E'_{n+3} = a' E_n + b' E_{n-1} + c' E_{n+1} \tag{94}$$

Proof. Informally we see that $E'_1 = a'$, $E'_2 = b'$, $E'_3 = c'$, $E'_4 = a' + c'$, $E'_5 = a' + b' + c'$, $E'_6 = a' + b' + 2c'$, $E'_7 = 2a' + b' + 3c'$ and so on.
Formally we prove the case *(n+1)* from the (known) cases *n* and *(n-2)*.

$E'_{n+4} = E'_{n+1} + E'_{n+3}$

$= a' E_{n-2} + b' E_{n-3} + c' E_{n-1} + a' E_n + b' E_{n-1} + c' E_{n+1}$

$= a'(E_{n-2} + E_n) + b'(E_{n-3} + E_{n-1}) + c'(E_{n-1} + E_{n+1})$

$= a' E_{n+1} + b' E_n + c' E_{n+2}$.

A similar proof works for the case $n \in \mathbb{Z}^-$ □

We denote the three given values from the sequence $(E'_n)_n$ as $E'_{n+3}, E'_{m+3},$ and $E'_{k+3}$ with n > m > k. From (94) we obtain the system of equations:

$$a'E_n + b'E_{n-1} + c'E_{n+1} = E'_{n+3}$$
$$a'E_m + b'E_{m-1} + c'E_{m+1} = E'_{m+3} \qquad (95)$$
$$a'E_k + b'E_{k-1} + c'E_{k+1} = E'_{k+3}$$

Using (79) we obtain the solutions (96):

$$a' = \frac{\begin{vmatrix} E'_{n+3} & E_{n-1} & E_{n+1} \\ E'_{m+3} & E_{m-1} & E_{m+1} \\ E'_{k+3} & E_{k-1} & E_{k+1} \end{vmatrix}}{\Delta_{n,m,k}}$$

$$b' = \frac{\begin{vmatrix} E_n & E'_{n+3} & E_{n+1} \\ E_m & E'_{m+3} & E_{m+1} \\ E_k & E'_{k+3} & E_{k+1} \end{vmatrix}}{\Delta_{n,m,k}}$$

$$c' = \frac{\begin{vmatrix} E_n & E_{n-1} & E'_{n+3} \\ E_m & E_{m-1} & E'_{m+3} \\ E_k & E_{k-1} & E'_{k+3} \end{vmatrix}}{\Delta_{n,m,k}}$$

with $\Delta_{n,m,k} = \Delta_{n,n+i,n+i+j} =: \Delta_{i,j} = E_i E_{i+j-1} - E_{i+j} E_{i-1}$, cf. (81). From this result we see that the necessary and sufficient condition for a solution is that $\Delta_{n,m,k} \neq 0$. From (94) we know that $\Delta_{n,m,k} = 0 \Leftrightarrow$ the fractions in (96) are of the form 0/0 and hence a', b', and c' cannot be determined. This happens e.g., for m=n+2 and k=m+1, m = n+5 and k=m+2, m = n+5 and k = m+3, and so on (cf. Table 5).

### 9.3.3 Examples

Let $(E'_n)_n$ be: 1,2,4,5,7,11,16,23,34,50,73, …[with the first value in this sequence, namely 1, corresponding with n = 1]

(i) Assume given $E'_8 = 23$, $E'_9 = 34$ and $E'_{10} = 50$. Then n+3 = 8 and hence n = 5 and i=j=1 so that $\Delta_{1,1} = 1$ (see (81)) and thus, by (96),

$$a' = \begin{vmatrix} 23 & 2 & 4 \\ 34 & 3 & 6 \\ 50 & 4 & 9 \end{vmatrix} = 1, b' = \begin{vmatrix} 3 & 23 & 4 \\ 4 & 34 & 6 \\ 6 & 50 & 9 \end{vmatrix} = 2; c' = \begin{vmatrix} 3 & 2 & 23 \\ 4 & 3 & 34 \\ 6 & 4 & 50 \end{vmatrix} = 4$$

These results correspond with the given $(E'_n)_n$ sequence.

(ii) Now we consider a case of non-consecutive indices. Assume given $E'_8 = 23$, $E'_9 = 34$ and $E'_{11} = 73$. Then n=5, i=1 and j=2. Now, $\Delta_{1,2} = 1$ and thus

$$a' = \begin{vmatrix} 23 & 2 & 4 \\ 34 & 3 & 6 \\ 73 & 6 & 13 \end{vmatrix} = 1, b' = \begin{vmatrix} 3 & 23 & 4 \\ 4 & 34 & 6 \\ 9 & 73 & 13 \end{vmatrix} = 2; c' = \begin{vmatrix} 3 & 2 & 23 \\ 4 & 3 & 34 \\ 9 & 6 & 73 \end{vmatrix} = 4$$

This again corresponds with the given $(E'_n)_n$ sequence.

(iii). Assume now that $E'_8 = 23$, $E'_{10} = 50$ and $E'_{11} = 73$. In this case, we have $E'_8 + E'_{10} = E'_{11}$, meaning that we only have two independent E' values. Here we have $\Delta_{i,j} = \Delta_{2,1} = 0$ and $a' = \frac{0}{0}$ as $\begin{vmatrix} 23 & 2 & 4 \\ 50 & 4 & 9 \\ 73 & 6 & 13 \end{vmatrix} = 0$.

Moreover also $\begin{vmatrix} 3 & 23 & 4 \\ 6 & 50 & 9 \\ 9 & 73 & 13 \end{vmatrix} = 0 = \begin{vmatrix} 3 & 2 & 23 \\ 6 & 4 & 50 \\ 9 & 6 & 73 \end{vmatrix}$ so that b' and c' too have the undetermined form $\frac{0}{0}$.

### 9.3.4 Consequences of formula (94)

Formula (94) yields the relation between an E'- sequence and a fixed E-sequence, if $a' = E'_1$, $b' = E'_2$, and $c' = E'_3$ are given. We will now use property (i)(a) and formula (82) which state that for every fixed i ∈ ℤ, $(\Delta_{i,j})_{j\in\mathbb{N}}$ is an E'-sequence with initial values given by (83), with b'= c'. Inserting these values in (94) we obtain, for all i ∈ ℤ fixed:

$$\Delta_{i,j+3} = a_{(i)}'E_j + b_{(i)}'(E_{j-1} + E_{j+1})$$
$$\Delta_{i,j+3} = a_{(i)}'E_j + b_{(i)}'E_{j+2} \tag{97}$$

This, together with (81) and (83) yields:

$$\Delta_{i,j+3} = (E_i^2 - E_{i-1}E_{i+1})E_j + (E_iE_{i+1} - E_{i-1}E_{i+2})E_{j+2}$$
$$= E_iE_{i+j+2} - E_{i+j+3}E_{i-1} \tag{98}$$

which is a separability property for $E_n$-values ($E_{i+j}$ as a function of $E_i$ and $E_j$).

### 9.3.5 Some observations about the case $\Delta_{i,j} = 0$

We note that (98) can be verified directly for j = 0,1,2, by using the classical equation $E_n + E_{n+2} = E_{n+3}$. We further see that this approach does not work when $\Delta_{i,j} = 0$. Yet, this happens only for a finite number of values $i,j \in \mathbb{N}$. This will be shown in the next theorem. Moreover, $E_n = 0$ only for a finite number of n, which are all smaller than or equal to zero. This theorem is important for its own sake.

**Theorem 17**

(a) The sequence $(E_n)_{n \in \mathbb{Z}}$ contains only a finite number of zeros, which all occur for n ≤ 0.

(b) The matrix $(\Delta_{i,j})_{i,j \in \mathbb{Z} \setminus \{0\}}$, with i+j ≠ 0, contains only a finite number of zeros.

Proof.

(a) We already know by (45) that $\forall n \in \mathbb{Z}: E_{-n-4} = E_{n+3}^2 - E_{n+2}E_{n+4}$
If we assume that $n \in \mathbb{N}$ then we also have:

$$E_{n-4} = E_{-n+3}^2 - E_{-n+2}E_{-n+4} \tag{99}$$

Applying now formula (17), which holds for all $n \in \mathbb{Z}$ we find:

$$E_{n-4} = C_1 a_0^{n-4} - C_1 \rho^{n-4} \cos((n-4)\zeta) + C_3 \rho^{n-4} \sin((n-4)\zeta)$$

$$E_{-n+2} = \frac{C_1}{a_0^{n-2}} - \frac{C_1}{\rho^{(n-2)}} \cos((n-2)\zeta) - \frac{C_3}{\rho^{(n-2)}} \sin((n-2)\zeta)$$

$$E_{-n+4} = \frac{C_1}{a_0^{n-4}} - \frac{C_1}{\rho^{(n-4)}} \cos((n-4)\zeta) - \frac{C_3}{\rho^{(n-4)}} \sin((n-4)\zeta)$$

where $a_0 \approx 1.4648493$, $\rho \approx 0.8265432$ and we investigate (99) for n tending to $+\infty$. We have:

$$\lim_{n \to +\infty} E^2_{-n+3} = \lim_{n \to +\infty} \left\{ C_1 a_0^{n-4} + \frac{X}{\rho^{2n-6}} \right\} \quad (100)$$

with $X = \bigl(C_1 \cos((n-2)\zeta) + C_3 \sin((n-2)\zeta)\bigr)\bigl(C_1 \cos((n-4)\zeta) + C_3 \sin((n-4)\zeta)\bigr)$. The expression X is clearly bounded.

Hence,

$$\lim_{n \to +\infty} E^2_{-n+3} = \lim_{n \to +\infty} a_0^{n-4} \left\{ C_1 + \frac{\left(\frac{1}{\rho^2}\right)^{n-3} X}{a_0^{n-4}} \right\}$$

As $\frac{1}{\rho^2} \approx 1.4637567 < a_0 \approx 1.4645493$, we find that

$$\lim_{n \to +\infty} E^2_{-n+3} = \lim_{n \to +\infty} C_1 a^{n-4} = +\infty \quad (100)$$

This result means that $E_{-n+3}$ for $n \in \mathbb{N}$ is bounded away from zero and hence $E_n$ $n \leq 0$ can have only a finite number of zeros. As $E_n > 0$ for $n \in \mathbb{N}$, this proves part (a) of this theorem.

(b) We already observed that we may take $i, j \in \mathbb{N}$, because $\forall\, i, j \in \mathbb{Z} \setminus \{0\}$, with $i+j \neq 0$, $\exists\, i', j' \in \mathbb{N}$ such that $\Delta_{i,j} = \pm \Delta_{i',j'}$. Indeed, $\Delta_{i,j} = \Delta_{n,m,k}$, with $n \neq m=n+i \neq k=m+j \neq n$ is independent of n and with $\Delta_{i',j'} = \Delta_{n',m',k'}$, where n'< m'< k', we have that $\Delta_{i',j'} = \pm \Delta_{i,j}$, dependent on the number of row changes needed to transform $\Delta_{n,m,k}$ into $\Delta_{n',m',k'}$. Hence it suffices to show that $\Delta_{i,j}= 0$ for a finite number of couples $(i,j) \in \mathbb{N} \times \mathbb{N}$.

We have, see (97), for $j \geq 1$ that

$$\Delta_{i,j+3}= a'_{(i)}E_j + b'_{(i)}E_{j+2}, \text{ with } a'_{(i)} = \Delta_{i,1}= E_{-i-1}, \; b'_{(i)} = \Delta_{i,2}= -E_{-i-3},$$

see (87),(88).

Hence,

$$\Delta_{i,j+3}= E_j \left( a'_{(i)} + b'_{(i)} \frac{E_{j+2}}{E_j} \right) \tag{101}$$

Now, by (100),

$$\lim_{i \to +\infty}(a'_{(i)})^2 = \lim_{i \to +\infty}(E_{-i-1})^2 = \lim_{i \to +\infty} (E_{-(i+4)+3})^2 = \lim_{i \to +\infty} C_1 a_0^{i}$$

Similarly, again using (100), we see :

$$\lim_{i \to +\infty}(b'_{(i)})^2 = \lim_{i \to +\infty}(E_{-i-3})^2 = \lim_{i \to +\infty}(E_{-(i+6)+3})^2 = \lim_{i \to +\infty} C_1 a_0^{i+2}$$

Hence: $\lim_{i \to +\infty}(\Delta_{i,j+3})^2 = C_1(E_j)^2 \lim_{i \to +\infty} \left( a_0^{i/2} \pm a_0^{(\frac{i}{2})+1} \frac{E_{j+2}}{E_j} \right)^2$

$$= C_1(E_j)^2 \left( \lim_{i \to +\infty} a^i \right) \left( 1 \pm a_0 \frac{E_{j+2}}{E_j} \right)^2$$

As $\forall j \geq 1$ (fixed), $a_0 \frac{E_{j+2}}{E_j} > 1$ the last factor in the expression above is not zero and hence $\lim_{i \to +\infty}(\Delta_{i,j+3})^2 = +\infty$ and thus $\Delta_{i,j} \neq 0$, $\forall i \geq i_0$, $\forall j \geq 4$ (fixed) and for a certain $i_0 \in \mathbb{N}$, independent from j. Now, for j =1,2,3, $\Delta_{i,1}= E_{-i-1}$, $\Delta_{i,2}= \Delta_{i,3}= -E_{-i-3}$, (cf. (87), (88), and hence also here there are a finite number of zeros (using part (a) of this theorem).

We still have to consider the case $\Delta_{i,j}$ for $1 \leq i \leq i_0$ and $j \in \mathbb{N}$. For fixed $i \in \{1, ..., i_0\}$ we know by (101) that

$$\lim_{j \to +\infty} \Delta_{i,j+3} = \lim_{j \to +\infty} E_j \left( a'_{(i)} + b'_{(i)} \frac{E_{j+2}}{E_j} \right) = \left( \lim_{j \to +\infty} E_j \right) \left( a'_{(i)} + b'_{(i)} a_0^2 \right)$$

by (18), with $a'_{(i)} = E_{-i-1}$, $b'_{(i)} = -E_{-i-3}$, independent from j. Now, $\forall i \in \{1, ..., i_0\}$, $a'_{(i)} + b'_{(i)} a_0^2 \neq 0$, because $a'_{(i)}, b'_{(i)} \in \mathbb{Z}$ and $a_0$ is irrational, it is impossible that $a_0^2 = -\frac{a'_{(i)}}{b'_{(i)}}$. We conclude that $\lim_{j \to +\infty} \Delta_{i,j+3} = +\infty$ or $-\infty$, implying that also in this case there are only a finite number of zeros.

### 9.4 The sequence $(D_n)_n$

We recall that the sequence $(D_n)_n \in \mathbb{Z}$ is the unique sequence determined by

$$D_n + D_{n+1} = D_{n+3} \tag{33}$$

with $D_{-1} = -1$, $D_0 = 1$, $D_1 = D_2 = 0$.

We denoted (equation (57)) the determinant

$$\begin{vmatrix} D_n & D_{n+2} & D_{n+1} \\ D_{n+1} & D_{n+3} & D_{n+2} \\ D_{n+2} & D_{n+4} & D_{n+3} \end{vmatrix}$$

as $\delta'_n$ and proved that, $\forall n \in \mathbb{Z}$, $\delta'_n = 1$ (58). This result too will be generalized considerably.

Definition

$$\Delta'_{n,m,k} =: \begin{vmatrix} D_{n+2} & D_{n+4} & D_{n+3} \\ D_{m+2} & D_{m+4} & D_{m+3} \\ D_{k+2} & D_{k+4} & D_{k+3} \end{vmatrix} \tag{100}$$

We see that for m = n+1, k=m+1, $\Delta'_{n,m,k} = \delta'_{n+2} = 1$. We denote in general m = n+i, k = m+j = n+i + j and $\Delta'_{n,m,k} = \Delta'_{i,j}$, again anticipating the fact that this determinant is independent of n.

Theorem 18. $\Delta'_{n,n+i,n+i+j}$ is independent of $n \in \mathbb{Z}$ and

$$\Delta'_{i,j} = D_{i+2}D_{i+j+4} - D_{i+j+2}D_{i+4} \tag{101}$$

Hence $\Delta'_{n,m,k}$ only depends on m-n=i and k-n=i+j.

Proof.

$$\Delta'_{n+1,n+i+1,n+i+j+1} = \begin{vmatrix} D_{n+3} & D_{n+5} & D_{n+4} \\ D_{n+i+3} & D_{n+i+5} & D_{n+i+4} \\ D_{n+i+j+3} & D_{n+i+j+5} & D_{n+i+j+4} \end{vmatrix}$$

Next, we apply the following equalities:

$$D_{n+5} = D_{n+2} + D_{n+3}$$
$$D_{n+i+5} = D_{n+i+2} + D_{n+i+3}$$
$$D_{n+i+j+5} = D_{n+i+j+2} + D_{n+i+j+3}$$

leading to:

$$\Delta'_{n+1,n+i+1,n+i+j+1} = \begin{vmatrix} D_{n+3} & D_{n+2} & D_{n+4} \\ D_{n+i+3} & D_{n+i+2} & D_{n+i+4} \\ D_{n+i+j+3} & D_{n+i+j+2} & D_{n+i+j+4} \end{vmatrix} = \Delta'_{n,n+i,n+i+j}$$

where we first exchanged the first and the second columns, and then the new second and third columns. From this result, we see that we may set n = 0 in (100). Hence,

$$\Delta'_{i,j} = \begin{vmatrix} 0 & 0 & 1 \\ D_{i+2} & D_{i+4} & D_{i+3} \\ D_{i+j+2} & D_{i+j+4} & D_{i+j+3} \end{vmatrix} = D_{i+2}D_{i+j+4} - D_{i+j+2}D_{i+4} \quad \square$$

We observe that for i=j=1 we find the known equality $\delta'_n = \Delta'_{1,1} = 1$. Properties similar to those proved for *E*- and *E'*-sequences, can here be shown too. We can further show the following application.

Problem. Given a D'-sequence $(D'_n)_n$ with initial values $a'' = D'_1$, $b'' = D'_2$, $c'' = D'_3$, determine *a''*, *b''* and *c''* from three given values in the sequence $(D'_n)_n$.

Solution

We start with a lemma.

Lemma 5. $\forall n \in \mathbb{Z}$:

$$D'_{n+3} = a'' D_{n+2} + b'' D_{n+4} + c'' D_{n+3} \tag{102}$$

Proof. Informally we see that: $D'_1 = a''$, $D'_2 = b''$, $D'_3 = c''$, $D'_4 = a'' + b''$, $D'_5 = b'' + c''$, $D'_6 = a'' + b'' + c''$, $D'_7 = a'' + 2b'' + c''$, $D'_8 = a'' + 2b'' + 2c''$ and so on.

Formally, proving the case for (n+1) from the known formula in the cases (n-1) and (n-2) we carry on as follows:

$$D'_{n+4} = D'_{n+1} + D'_{n+2}$$
$$= a''D_n + b''D_{n+2} + c''D_{n+1} + a''D_{n+1} + b''D_{n+3} + c''D_{n+2}$$
$$= a''(D_n + D_{n+1}) + b''(D_{n+2} + D_{n+3}) + c''(D_{n+1} + D_{n+2})$$
$$= a''D_{n+3} + b''D_{n+5} + c''D_{n+4}$$

A similar proof works for $n \in \mathbb{Z}^-$. □

Denoting the given values in the $(D'_n)_n$ − sequence as $D'_{n+3}, D'_{m+3}, D'_{k+3}$ with n > m > k we obtain from (102):

$$a''D_{n+2} + b''D_{n+4} + c''D_{n+3} = D'_{n+3}$$
$$a''D_{m+2} + b''D_{m+4} + c''D_{m+3} = D'_{m+3} \tag{103}$$

$$a"D_{k+2} + b"D_{k+4} + c"D_{k+3} = D'_{k+3}$$

Hence, if $\Delta'_{n,m,k} \neq 0$, we find:

$$a" = \frac{\begin{vmatrix} D'_{n+3} & D_{n+4} & D_{n+3} \\ D'_{m+3} & D_{m+4} & D_{m+3} \\ D'_{k+3} & D_{k+4} & D_{k+3} \end{vmatrix}}{\Delta'_{n,m,k}}, b" = \frac{\begin{vmatrix} D_{n+2} & D'_{n+3} & D_{n+3} \\ D_{m+2} & D'_{m+3} & D_{m+3} \\ D_{k+2} & D'_{k+3} & D_{k+3} \end{vmatrix}}{\Delta'_{n,m,k}}, c = \frac{\begin{vmatrix} D_{n+2} & D_{n+4} & D'_{n+3} \\ D_{m+2} & D_{m+4} & D'_{m+3} \\ D_{k+2} & D_{k+4} & D'_{k+3} \end{vmatrix}}{\Delta'_{n,m,k}}$$

(104)

Similar examples as in the E'-case can also here be given (left to the reader). We conclude with the following simple result, similar to the result obtained in the F'-case:

Theorem. The sequence $(E'_n)_n$ is the right-hand side of the $(E_n)_n$ – sequence if and only if there exists $n_0 \in \mathbb{Z}$ such that $a' =: E'_1 = E_{n_0}$, $b' =: E'_2 = E_{n_0+1}$ ; $c' =: E'_3 = E_{n_0+2}$. A similar result holds for the sequence $(D'_n)_n$.

We note that it is easy to see that, except for the zero-sequence no $(E'_n)_n$ – sequence, can be a $(D'_n)_n$-sequence.

We did not find names associated with $(E'_n)_n$ – sequences, but for $(D'_n)_n$- sequences we found the names Perrin for the case $D'_1 = a" = 0$, $D'_2 = b" = 2; D'_3 = c" = 3$ and Padovan for the case $D'_1 = D'_2 = 1, D'_3 = 2$.

## 10. Conclusion

In the first part of this study, we showed that for a continuous function f of positive variables the requirement $h_f = f$ leads to the conclusion that either f = **0**, or $f(x) = x^{\left(\frac{\sqrt{5}+1}{2}\right)}$, leading to a relation with the golden ratio. A similar requirement for discrete values leads to f = {(1,1)} and no other solution. The reason for these results is that $h_f$ = f implies that f(x)/x must belong to the domain of f. This is not a restriction on $\mathbb{R}_0^+$, but it is a severe one on a set of discrete values. Yet, in the two investigations, Fibonacci numbers played a role, resulting from the fact that $h_f$ = f implies that f(f(x)) = x f(x).

We also studied the relation $h_f = f \circ f$ and found the same results as in the $h_f = f$ case, both for continuous f as well as in the discrete case but the proofs are considerably more difficult than in the Fibonacci case. Now the Narayana cows sequence $(E_n)_{n \in \mathbb{Z}}$ is involved, featuring the characteristic determinant

$$\Delta_{i,j} := \begin{vmatrix} E_n & E_{n-1} & E_{n+1} \\ E_{n+i} & E_{n+i-1} & E_{n+i+1} \\ E_{n+i+j} & E_{n+i+j-1} & E_{n+i+j+1} \end{vmatrix}$$

It was shown that $\Delta_{i,j}$ is independent of $n \in \mathbb{Z}$ and equals

$$E_i E_{i+j-1} - E_{i+j} E_{i-1}$$

from which several identities involving $(E_n)_{n \in \mathbb{Z}}$ were proved.

This determinant also plays a key role in determining generalized Narayana sequences, based on the classical $(E_n)_{n \in \mathbb{Z}}$ sequence.

We further showed that $E_n$ = 0 and $\Delta_{i,j}$ = 0 for only a finite number of integers *n* and pairs *(i,j)* $\in \mathbb{Z} \times \mathbb{Z}$.

Similarly, we studied the relation $h_{f \circ f} = f$ resulting in a Padovan-type sequence $(D_n)_{n \in \mathbb{Z}}$. Now the key determinant is

$$\Delta'_{i,j} = \begin{vmatrix} D_n & D_{n+2} & D_{n+1} \\ D_{n+i} & D_{n+i+2} & D_{n+i+1} \\ D_{n+i+j} & D_{n+i+j+2} & D_{n+i+j+1} \end{vmatrix}$$

Also $\Delta'_{i,j}$ is independent of $n \in \mathbb{Z}$ and equals

$$D_{i+2}D_{i+j+4} - D_{i+j+2}D_{i+4}$$

Similar properties and applications as in the $(E_n)_{n \in \mathbb{Z}}$ case were proved.

We end this work with some open problems.

Open problem 1. Prove (67) without the goniometric forms of $E_n$ and $D_n$.

Open problem 2. Find explicit formulae to describe the relationship between $D_n$, $E_n$, and $F_n$ ($n \in \mathbb{Z}$).

Open problem 3. Find those $i, j \in \mathbb{N}$ or $\mathbb{Z}$ such that $\Delta_{i,j} = 0$, and the same problem for $\Delta'_{i,j} = 0$.

**Acknowledgment.** The author thanks Ronald Rousseau for his advice in the solution of two difference equations and for editing this article. He also thanks Li Li (National Science Library, CAS) for drawing Figure 1.